\PassOptionsToPackage{hyphens}{url}
\documentclass[journal]{IEEEtran}

\newcommand{\beq}{\begin{equation}}
\newcommand{\eeq}{\end{equation}}
\newcommand{\bsq}{\begin{subequations}}
\newcommand{\esq}{\end{subequations}}
\newcommand{\bq}{\begin{eqnarray}}
\newcommand{\eq}{\end{eqnarray}}
\newcommand{\bqn}{\begin{eqnarray*}}
\newcommand{\eqn}{\end{eqnarray*}}
\usepackage{bbm}
\usepackage[pdftex]{graphicx}
\graphicspath{{Figs/}}
\usepackage{enumerate}
\usepackage{enumitem} 

\usepackage{mathptmx}       
\DeclareMathAlphabet{\mathcal}{OMS}{cmsy}{m}{n}
\usepackage{helvet}         
\usepackage{courier}        
\usepackage{type1cm}        
\usepackage{makeidx}         
\usepackage{graphicx}        
\usepackage{url} 
\usepackage{multicol}        
\usepackage[bottom]{footmisc}
\usepackage{ifthen}
\usepackage{subfigure}
\usepackage[usenames,dvipsnames]{color}

\makeindex             

\usepackage{graphics} 
\usepackage{graphicx} 
\usepackage{epsfig} 
\usepackage{epstopdf}
\usepackage{mathptmx} 
\usepackage{times} 
\usepackage[cmex10]{amsmath}

\usepackage{mathtools}  
\usepackage{amssymb}  
\usepackage{amsthm}

\usepackage{amsfonts}
\usepackage{cite}
\usepackage{verbatim}
\usepackage{comment}
\usepackage{algorithmic}
\usepackage{multirow}
\usepackage{cite}
\usepackage{stfloats}
\usepackage{supertabular}
\usepackage{longtable}

\usepackage[hyphens]{url}
\usepackage{breakurl}

\DeclareGraphicsExtensions{.pdf,.png,.jpg,.eps,.ps}
\renewcommand{\arraystretch}{1.2}
\usepackage{arydshln}
\usepackage{multicol}
\usepackage{colortbl}
\theoremstyle{definition}

\newtheorem{proposition}{Proposition}

\theoremstyle{definition}
\newtheorem{definition}{Definition}

\ifCLASSINFOpdf
\else
\fi

\usepackage{comment}
\usepackage{ifthen}
\newboolean{showcomments}
\setboolean{showcomments}{true}
\newcommand{\ychen}[1]{\ifthenelse{\boolean{showcomments}}
        { \textcolor{red}{YC: #1}}}
\newcommand{\tongxin}[1]{\ifthenelse{\boolean{showcomments}}
        { \textcolor{black}{(#1)}}{}}

\usepackage{amsmath}
\allowdisplaybreaks[4]

\usepackage[ruled]{algorithm2e}
\usepackage{booktabs}

\begin{document}

%
\title{Robust Microgrid Dispatch with Real-Time Energy Sharing and Endogenous Uncertainty}

\author{
Meng~Yang,
Rui~Xie,~\IEEEmembership{Member,~IEEE},
Yongjun Zhang,~\IEEEmembership{Senior Member,~IEEE},
Yue~Chen,~\IEEEmembership{Member,~IEEE}
\thanks{This work was supported by the National Natural Science Foundation of China under Grant No. 52307144. (Corresponding to Y. Chen)}
\thanks{M. Yang, R. Xie, and Y. Chen are with the Department of Mechanical and Automation Engineering, the Chinese University of Hong Kong, Hong Kong SAR, China (e-mail: myang@mae.cuhk.edu.hk, ruixie@cuhk.edu.hk, yuechen@mae.cuhk.edu.hk).}
\thanks{Y. Zhang is with the School of Electric Power Engineering, South China University of Technology, Guangzhou, China. (email: zhangjun@scut.edu.cn)}
}

\markboth{Journal of \LaTeX\ Class Files,~Vol.~XX, No.~X, Feb.~2019}%
{Shell \MakeLowercase{\textit{et al.}}: Bare Demo of IEEEtran.cls for IEEE Journals}
%



\maketitle

\begin{abstract}
With the rising adoption of distributed energy resources (DERs), microgrid dispatch is facing new challenges: DER owners are independent stakeholders seeking to maximize their individual profits rather than being controlled centrally; and the dispatch of renewable generators may affect the microgrid’s exposure to uncertainty. To address these challenges, this paper proposes a two-stage robust microgrid dispatch model with real-time energy sharing and endogenous uncertainty. In the day-ahead stage, the connection/disconnection of renewable generators is optimized, which influences the size and dimension of the uncertainty set. As a result, the uncertainty set is endogenously given. In addition, non-anticipative operational bounds for energy storage (ES) are derived to enable the online operation of ES in real-time.  In the real-time stage, DER owners (consumers and prosumers) share energy with each other via a proposed energy sharing mechanism, which forms a generalized Nash game. To solve the robust microgrid dispatch model, we develop an equivalent optimization model to compute the real-time energy sharing equilibrium. Based on this, a projection-based column-and-constraint generation (C\&CG) method is proposed to handle the endogenous uncertainty. Numerical experiments show the effectiveness and advantages of the proposed model and method.
\end{abstract}


%
\IEEEpeerreviewmaketitle

\section*{Nomenclature}
\addcontentsline{toc}{section}{Nomenclature}
\subsection{Abbreviation}
\begin{IEEEdescription}[\IEEEusemathlabelsep\IEEEsetlabelwidth{ssssssss}]
\item[C\&CG] Column-and-constraint generation.
\item[DER] Distributed energy resource.
\item[RMD] Robust microgrid dispatch.
\item[RG] Renewable generator.
\item[GNE] Generalized Nash equilibrium.
\item[SOC] State-of-charge.
\end{IEEEdescription}
\subsection{Indices and Sets}
\begin{IEEEdescription}[\IEEEusemathlabelsep\IEEEsetlabelwidth{ssssssss}]
\item[$k \in \mathcal{I}$] Set of consumers.
\item[$k \in \mathcal{J}$] Set of prosumers.
\item[$k \in \mathcal{K}$] Set of customers.
\item[$t \in \mathcal{T}$] Set of periods.
\item[$g \in \mathcal{G}$] Set of dispatchable gas-fired units.
\item[$e \in \mathcal{E}$] Set of energy storage units.
\item[$b \in \mathcal{B}$] Set of nodes.
\item[$(n,b) \in \mathcal{L}$] Set of lines.
\end{IEEEdescription}
       
\subsection{Parameters}
\begin{IEEEdescription}[\IEEEusemathlabelsep\IEEEsetlabelwidth{ssssssss}]
\item[$\underline{W}_{kt},\overline{W}_{kt}$] Lower/upper bound of the real-time output of RG $k$ in period $t$.
\item[$\underline{P}_{g},\overline{P}_{g}$] Lower/upper bound of the available active power of dispatchable gas-fired units $g$.
\item[$\underline{Q}_{g},\overline{Q}_{g}$] Lower/upper bound of the available reactive power of dispatchable gas-fired units $g$.
\item[$\eta_{e}^{c},\eta_{e}^{d}$] Charging/discharging efficiency of ES $e$.
\item[${d}_{kt}^{f}$] Fixed demand of customer $k$.
\item[$\underline{D}_{kt},\overline{D}_{kt}$] Lower/upper bound of customer $k$'s elastic demand in period $t$.
\item[$\overline{P}_{nb},\overline{Q}_{nb}$] Upper bounds of the active/reactive power of transmission line $(n,b)$.
\item[$U_k(.)$] Disutility function of customer $k$.
\item[$a$] Market sensitivity.
\item[$r_{nb},x_{nb}$] Resistance/reactance of line $(n,b)$.
\item[$P_{bt}^{0},Q_{bt}^{0}$] Active/reactive power load of node $b$ in period $t$.
\item[$B^S$] Uncertainty budget among different RGs.
\item[$B^T$] Uncertainty budget among different periods.
\item[$\alpha_{kt}$] Unit penalty cost of underutilized renewable energy $k$ in period $t$.
\item[$c_g$] Unit operation cost of dispatchable gas-fired units including fuel cost and carbon intensity cost.
\item[$s_g$] Unit reserve cost of dispatchable gas-fired units.
\end{IEEEdescription}

\subsection{Decision Variables}
\begin{IEEEdescription}[\IEEEusemathlabelsep\IEEEsetlabelwidth{ssssssssss}]
\item[$u_{kt}$] Binary decision variable to indicate if RG $k$ is connected in period $t$.
\item[$h_{gt}$] Binary decision variable to indicate if dispatchable gas-fired units $g$ is on in period $t$.
\item[$\mu_{et}^{c},\mu_{et}^{d}$] Binary decision variable to indicate if ES $e$ is charged/discharged in period $t$.
\item[$\underline{E}_{et}^{DA},\overline{E}_{et}^{DA}$] SOC bounds of ESs.
\item[$\underline{p}_{et}^{DA,c},\overline{p}_{et}^{DA,c}$] Charging power bounds of ESs.
\item[$\underline{p}_{et}^{DA,d},\overline{p}_{et}^{DA,d}$] Discharging power bounds of ESs.
\item[${p}_{et}^{c},{p}_{et}^{d}$] Real-time charging/discharging power of ESs.
\item[$w_{kt}$] Real-time output of RG $k$ in period $t$.
\item[$P_{gt},Q_{gt}$] Day-ahead active/reactive output of dispatchable gas-fired units $g$ in period $t$.
\item[$\Delta_{gt},r_{gt}$] Real-time adjustment output and operational reserve of dispatchable gas-fired units $g$ in period $t$.
\item[$q_{kt}$] Sharing quantity of customer $k$ in period $t$.
\item[$\lambda_{kt}$] Sharing price of customer $k$ in period $t$.
\item[$b_{kt}$] Bid of customer $k$ in the energy sharing market.
\item[$v_{bt}$] Voltage of node $b$ in period $t$.
\end{IEEEdescription}

\section{Introduction}

%
%
%
%

\IEEEPARstart{T}{he} rapid growth of distributed energy resources (DERs) brings significant changes to power systems. Though DERs can reduce the dependence of power systems on fossil fuel and alleviate environmental pollution, their fluctuating nature exacerbates the difficulty of maintaining real-time energy balancing \cite{ellabban2014renewable}. Meanwhile, consumers with DERs are turning into prosumers that can provide more flexibility to power systems by participating in energy management proactively \cite{parag2016electricity}. Therefore, there is an urgent call for effective methods for managing active prosumers at demand side to accommodate volatile renewable energy. Robust microgrid dispatch (RMD) is a widely used approach \cite{qiu2020recourse} for handling renewable uncertainties in microgrids. However, traditional RMD models face two major challenges due to the proliferation of DERs.

The first challenge is that while traditional RMD models concentrated on exogenous uncertainty independent of the first-stage decision variables \cite{qiu2020tri,chu2021frequency}, endogenous uncertainty is becoming more common under DER integration. Particularly, with the increasingly frequent oversupply of renewable generation in countries such as Poland, America, and Australia, distributed renewable generators are required to be capable of being remotely disconnected to ensure power system security \cite{olivella2018optimization,Website}. The connection/disconnection of distributed renewable generators is an important decision that leads to endogenous uncertainty \cite{chen2022robust}, which is considered in this paper. Endogenous uncertainty has attracted great attention in recent years and a thorough review on related studies was provided in \cite{wang2024application}. Endogenous uncertainty was first considered under a stochastic programming framework. Reference \cite{qi2023chance} addressed the chance-constrained energy storage operation problem under endogenous uncertainty by an iterative algorithm. An affine function-based solution method \cite{yin2022coordinated} and a modeling transformation technique \cite{qi2023chance} were proposed to solve the stochastic model with endogenous uncertainty. Since security is of top priority in power system operation, addressing endogenous uncertainty in a robust optimization (RO) framework is also crucial \cite{zeng2022two}. Traditional RO algorithms such as column-and-constraint generation (C\&CG) \cite{zeng2013solving} may not be applicable to RO with endogenous uncertainty, since the previously selected worst-case scenarios may fall outside of the updated uncertainty set. Reformulation approaches to solving static RO with endogenous uncertainty have been developed \cite{nohadani2018optimization}-\cite{lappas2018robust}. Furthermore, adjustable robust optimization (ARO) problems that accounts for the adjustment after uncertainty realizations have been the recent focus. The solution to ARO with endogenous uncertainty can be estimated via a K-adaptability approximation method \cite{vayanos2020robust}. Modified Benders decomposition method \cite{zhang2021robust}, multi-parametric programming \cite{avraamidou2020adjustable}, adaptive C\&CG algorithm \cite{su2022multi}, dual C\&CG algorithm \cite{tan2024robust}, improved column generation algorithm \cite{chen2023robust}, and a parametric variant of C\&CG algorithm \cite{zeng2022two} were developed to offer an exact solution. While the above studies have offered valuable methods to solve RO with endogenous uncertainty, they mostly focused on continuous uncertainty set and can hardly handle the uncertainty set with binary variables. Furthermore, their lower-level problems are centralized optimization models while the one in this paper is a generalized Nash game. The unique features of the problem studied require a novel design of the solution algorithm. 

Traditional RMD models also face the other critical challenge: Different from conventional controllable units, prosumers at the demand side, such as smart buildings \cite{lyu2021fully}, smart homes \cite{celik2017decentralized}, and electric vehicle charging stations \cite{lv2023optimal} are self-interested stakeholders, which can dispatch DERs and participate in energy trading to maximize their own profits. The centralized dispatch scheme may become inapplicable due to the competing interests between prosumers and the operator. Innovative market mechanisms are urgently needed to coordinate prosumers to act in a way that achieves the operator's objective. Peer-to-peer energy sharing as a promising remedy has attracted widespread attention in recent years \cite{chen2022review}. Many countries have established demonstration projects, such as the GridWise project funded by the U.S. Department of Energy \cite{forfia2016view}, the Couperus project in the suburbs of The Hague, the Netherlands \cite{gjorgievski2021potential} and the distributed generation market trading pilot project in Zhenglu Industrial Park, Changzhou City, China. Reference \cite{tushar2021peer} provided a detailed summary of pilot projects on P2P energy sharing in North America, Europe, Australia, and Asia. Studies on energy sharing mechanism design can be divided into game theoretic-based ones and optimization-based ones. Shapley value\cite{mei2019coalitional}, nucleolus\cite{yang2021optimal}, and Nash bargaining\cite{cui2021economic} are three common methods to allocate the energy sharing costs/benefits in a cooperative game setting. Noncooperative game-based mechanisms using Stackelberg game\cite{xu2020data}, generalized Nash game\cite{chen2018analyzing}, multi-leader multi-follower game\cite{anoh2019energy}, and bilateral Nash game\cite{morstyn2018bilateral} were also developed. Besides the game theoretic-based mechanisms, which may result in a suboptimal equilibrium, the optimization-based mechanisms can achieve social optimum via distributed optimization frameworks \cite{wang2019incentivizing}. However, few of existing research has integrated energy sharing into a two-stage decision-making framework and quantified its ability to accommodate volatile renewable energy. Reference \cite{xu2019distributed} proposed a distributed robust energy management system for microgrids, where each microgrid solves a robust scheduling problem and the energy sharing among networked microgrids is settled via alternating direction method of multipliers (ADMM) algorithm. The energy sharing between prosumers served as the lower level of a robust model in \cite{cui2018two}, considering the uncertainties from renewable energy and market prices. The two-settlement transactive energy sharing mechanism for DER aggregators was developed in \cite{wang2022transactive}, adopting ARO to model the local problem with multi-source uncertainties. However, the above works focused on exogenous uncertainty. Moreover, the energy sharing between prosumers was formulated as a centralized optimization problem, which can be improved by using a game model instead to better characterize prosumers' interactions.

In this paper, we propose a novel RMD model and its solution algorithm to deal with the above two challenges, i.e. \emph{endogenous uncertainty} due to renewable connection/disconnection in the day-ahead stage and \emph{real-time energy sharing} to accommodate proactive prosumers. Our main contributions are two-fold:

1) \textbf{Energy sharing embedded RMD model with endogenous uncertainty}. We propose a RMD model with endogenous uncertainty and real-time energy sharing. In the day-ahead stage, the microgrid operator decides on the connection/disconnection strategies of distributed renewable generators and the lower and upper bounds of charging/discharging power and state-of-charge (SOC) for energy storage (ES) units. The connection/disconnection strategies influence the range of real-time renewable power output variations, leading to an endogenous uncertainty set; the ES operational bounds ensure the feasibility of non-anticipative real-time ES operation. In the real-time stage, customers (consumers and prosumers) are allowed to share energy with each other to deal with the renewable output deviations. An energy sharing mechanism is proposed to facilitate this energy exchange, under which all customers play a generalized Nash game.

2) \textbf{Solution algorithm based on equilibrium characterization and projection-based C\&CG}. We prove the existence of the energy-sharing market equilibrium and offer a convex centralized optimization to solve it. Based on this, the RMD model with a real-time energy sharing game turns into an adjustable robust optimization with endogenous uncertainty, where the conventional C\&CG algorithm gets stuck. Based on the structure of the endogenous uncertainty set, a projection-based C\&CG algorithm is developed to solve the problem efficiently. The proposed algorithm can consider customers' self-interested strategies while enabling the operator to solve the RMD problem centrally. Case studies demonstrate the efficiency and scalability of the proposed algorithm. Some interesting phenomena are also revealed.

The rest of this paper is organized as follows: The real-time energy sharing mechanism and the two-stage RMD model are presented in Section \ref{sec-2}. In Section \ref{sec-3}, an equivalent optimization model is provided to calculate the energy sharing equilibrium, and a projection-based C\&CG algorithm is developed to handle the endogenous uncertainty. Case studies are given in Section \ref{sec-4} with conclusions in Section \ref{sec-5}.

\section{Mathematical Models}
\label{sec-2}
\begin{figure}[t]
\centering
\includegraphics[width=0.6\columnwidth]{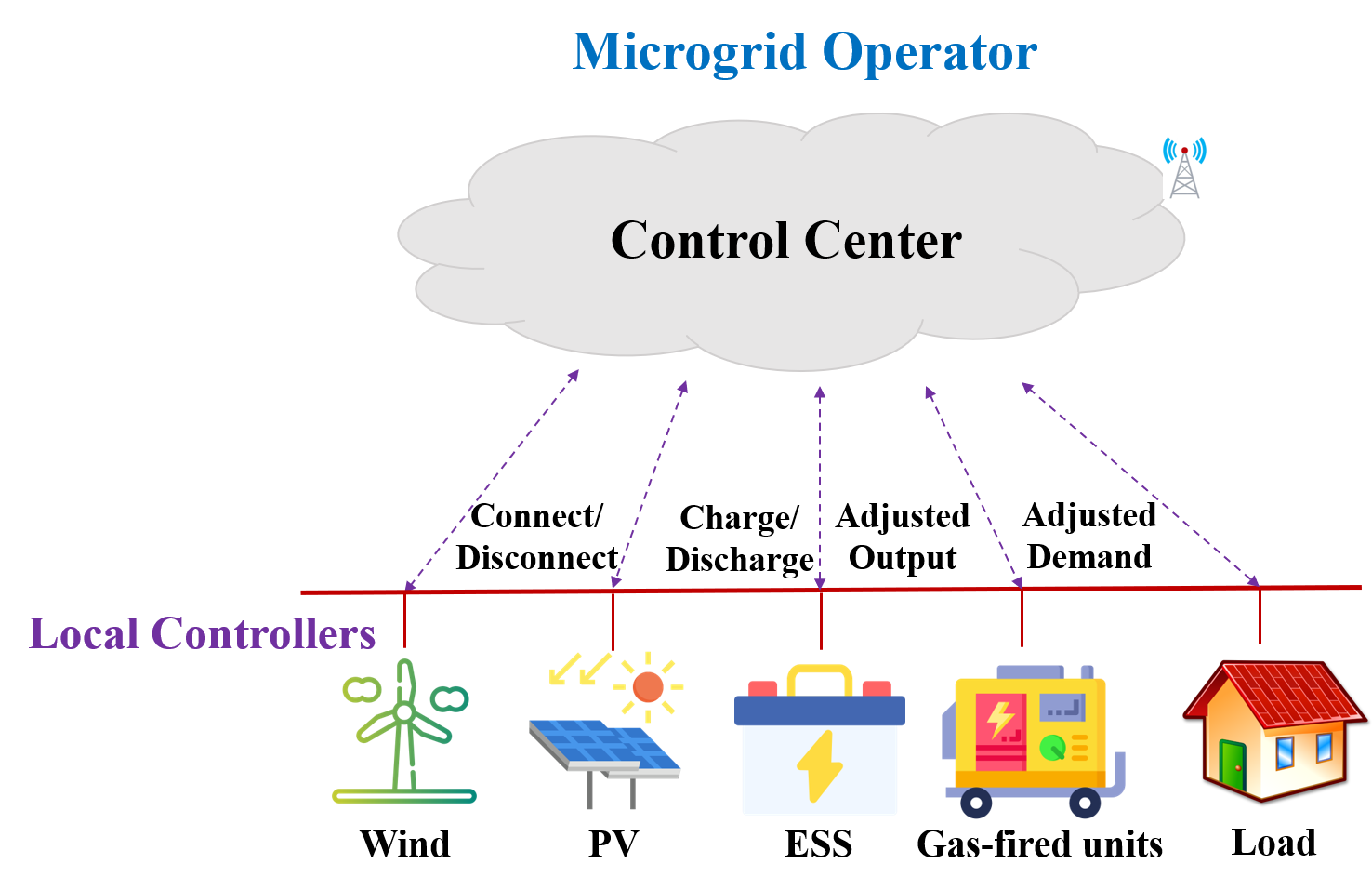}
\caption{Overview of the control architecture in a microgrid.}
\label{fig:DER-Controller}
\end{figure}

We consider the dispatch of a standalone microgrid with controllable gas-fired units, energy storage, consumers, and prosumers. Each consumer indexed by $k \in \mathcal{I}=\{1,...,I\}$ has load only, while each prosumer indexed by $k \in \mathcal{J}=\{1,...,J\}$ has both load and a disconnectable RG. Some loads are fixed while the others are adjustable. Denote by $\mathcal{K}$ the union of sets $\mathcal{I}$ and $\mathcal{J}$, i.e., $\mathcal{K}:=\mathcal{I} \cup \mathcal{J}$. As illustrated in Fig. \ref{fig:DER-Controller}, the microgrid operator can send control commands including connection/disconnection of each distributed RG, charged/discharged status of each energy storage unit, and power setpoints of each controllable gas-fired unit and load to the local controllers. This setting is practical thanks to the rapid development of intelligent electronic devices and information and communication technology \cite{ross2015multiobjective}. A similar setting was adopted in references \cite{qi2016cybersecurity,olivella2018optimization}.


\begin{figure}[t]
\centering
\includegraphics[width=1.0\columnwidth]{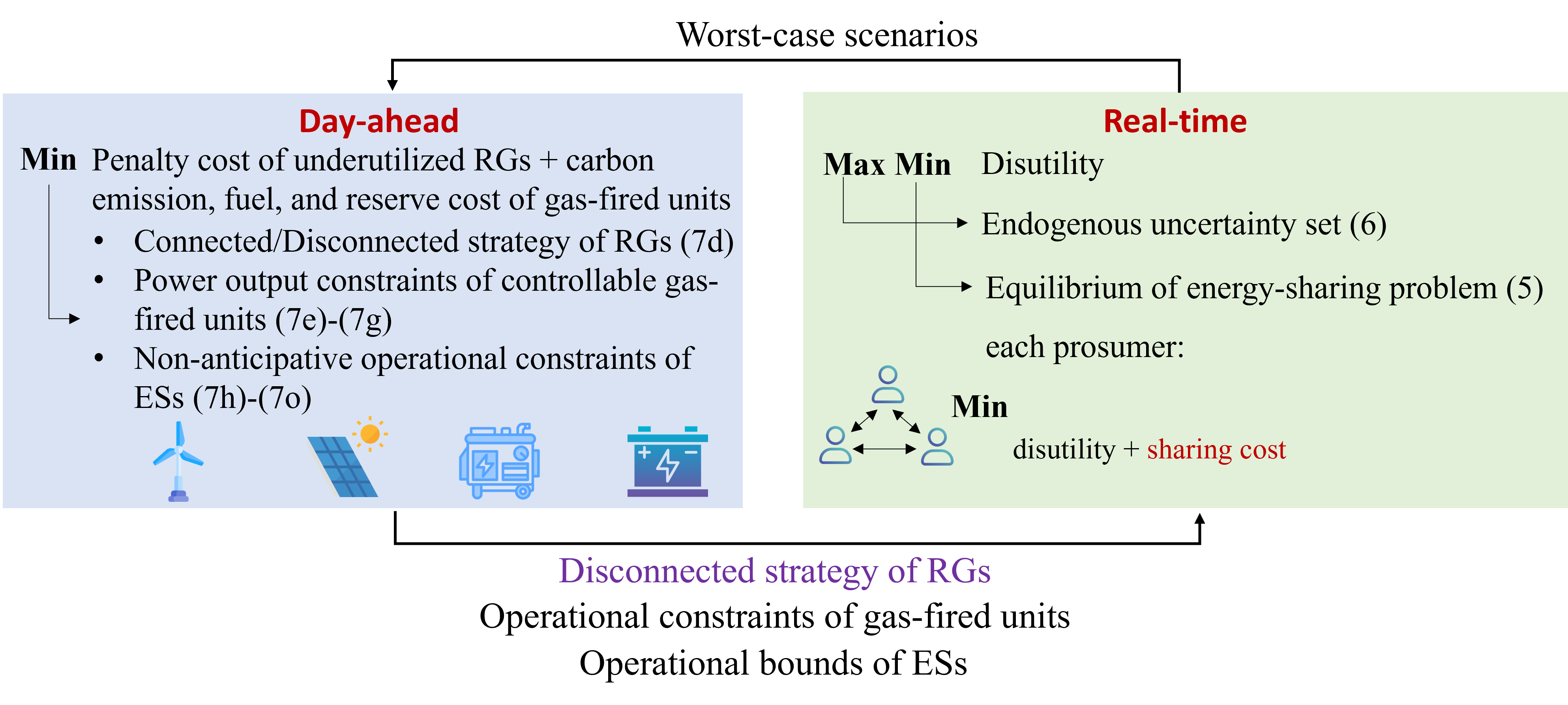}
\caption{The proposed two-stage decision-making framework.}
\label{fig:structure}
\end{figure}

In the day-ahead stage, the microgrid operator decides on the connection/disconnection of RGs for $T$ hours. The power outputs of connected RGs are uncertain, varying within certain ranges. In the real-time stage, to deal with RG power deviations, customers can adjust their power or share energy with each other. In addition, the outputs of controllable gas-fired units and energy storage units can be adjusted as well. Considering that the microgrid operates in an online manner period-by-period in real-time, it is difficult to dispatch the time-coupled energy storage units efficiently. To address this challenge, we propose to set non-anticipative operational bounds for energy storage units in the day-ahead stage. This enables time-decoupled operation of energy storage units in real-time by adjusting their charging and discharging power within the operational bounds for each period. The two-stage decision-making framework is given in Fig. \ref{fig:structure}.  In the following, we first model the real-time stage, and then the endogenous uncertainty set, and finally the two-stage RMD model.


\subsection{Real-Time Stage with Energy Sharing}
\label{sec-2-1}
In each period $t \in \mathcal{T}=\{1,...,T\}$, suppose the real-time RG output is $w_{kt}$ for unit $ k \in \mathcal{J}$. $u_{kt}$ is the day-ahead binary decision indicating whether the RG $k$ is connected or not. If $u_{kt}=0$, we have $w_{kt}=0$; otherwise, the real-time RG output is a random value within $[\underline{W}_{kt}, \overline{W}_{kt}]$. The customer can alter its demand $d_{kt}$ within $[\underline{D}_{kt},\overline{D}_{kt}]$
and exchange $q_{kt}$ with other customers at a price $\lambda_{kt}$ to help balance the real-time power. Each customer tries to minimize its individual cost, i.e. its payment to the energy sharing market $\lambda_{kt}q_{kt}$ plus its disutility due to the demand adjustment. The disutility can be quantified by a strictly convex function $U_{k}(d_{kt})$.

Here, we propose an energy sharing mechanism considering network constraints. Let $b_{kt}$ be the bid of customer $k \in \mathcal{K}$ at period $t \in \mathcal{T}$. A generalized demand function \cite{hobbs2000strategic}-\cite{li2015demand} is used to describe the relationship between $q_{kt}$, $\lambda_{kt}$, and $b_{kt}$:
\begin{align}\label{eq:df}
    q_{kt}~=~ -a\lambda_{kt}+b_{kt},\forall k \in \mathcal{K}, \forall t \in \mathcal{T},
\end{align}
where $a>0$ is the market sensitivity, and $q_{kt}>0$ indicates a buyer while $q_{kt}<0$ indicates a seller. 

\emph{Remark on the demand function:} In practice, the market sensitivity $a$ can be estimated by machine learning methods using historical demand and price data. We prove later in Proposition \ref{Thm:prop-central} that even when the parameter $a$ does not match each prosumer’s individual price sensitivity accurately, the desired properties of the market (e.g. existence and uniqueness of the equilibrium) still exist. In fact, $a$ can be regarded as a parameter in the set rule for sharing market clearing. Moreover, a linear demand function \eqref{eq:df} is used in this paper. A nonlinear demand function may exist in some highly dynamic market environment. In that case, we may iteratively linearize it using the first-order Taylor expansion near the updated price.

For each period $t \in \mathcal{T}$, the energy sharing market works as follows:

\textbf{Step 1:} (Initialization) Estimate the market sensitivity $a$ via historical data. Each customer inputs its private information $U_{k}(.)$, $\underline{D}_{kt}$, $\overline{D}_{kt}$ to its smart meter. Set $\lambda_{kt}^1=0,\forall k \in \mathcal{K}$, and $n=1$. Choose tolerance $\epsilon$.

\textbf{Step 2:} Based on the latest sharing price $\lambda_{kt}^n$, each customer $k \in \mathcal{K}$ determines the optimal bid $b_{kt}^{n+1}$ by solving  \eqref{eq:sharing-pro}, and submits the bid to the sharing platform.
\bsq
\label{eq:sharing-pro}
\begin{align}
    \mathop{\min}_{d_{kt},b_{kt},q_{kt},\forall t \in \mathcal{T}}~ & \Gamma_k^c(d_{kt},q_{kt},\lambda_{kt}) \label{eq:sharing-pro.1}\\
    \mbox{s.t.} ~&  q_{kt}=-a\lambda_{kt}+b_{kt} ,\label{eq:sharing-pro.2}\\
    ~ & \left\{\begin{aligned}  
     & q_{kt}=d_{kt}+d_{kt}^{f},~\mbox{if}~ k \in \mathcal{I} \\
     & w_{kt}+q_{kt}=d_{kt}+d_{kt}^{f},~\mbox{if}~ k \in \mathcal{J}
    \end{aligned}
    \right.,
    \label{eq:sharing-pro.3}\\
    ~ & \underline{D}_{kt} \le d_{kt} \le \overline{D}_{kt},\label{eq:sharing-pro.4}
\end{align}
\esq
where $\Gamma_k^c(d_{kt},q_{kt},\lambda_{kt}):= \left[U_k(d_{kt})+ \lambda_{kt}q_{kt}\right]$ is the total cost. Constraint \eqref{eq:sharing-pro.3} is the power balance condition and \eqref{eq:sharing-pro.4} limits the adjustable range of demand. Denote the feasible set of problem \eqref{eq:sharing-pro} as $\mathcal{X}_{kt}^c(\lambda_{kt})$.

\textbf{Step 3:} Upon receiving all the bids $b_{kt}^{n+1},\forall k \in \mathcal{K}$, the microgrid operator updates the energy sharing prices $\lambda_{kt}^{n+1},\forall k \in \mathcal{K}$ by solving: 
\bsq
\label{eq:sharing-oper}
\begin{align}
    \mathop{\min}_{\lambda_{kt},\forall k}~ & \sum \nolimits_{k \in \mathcal{K}} \lambda_{kt}^2 ,\label{eq:sharing-oper.1}\\
    \mbox{s.t.}~ & (-a\lambda_{kt}+b_{kt}^{n+1}) \in \mathcal{F}_{c}, \label{eq:sharing-oper.2}
\end{align}
\esq
where 
\bsq
\label{eq:Fc}
\begin{align}
\mathcal{F}_c & := \{ q_{kt}, \forall k ~|~ \exists p_{et}^c, p_{et}^d, P_{bt}, Q_{bt}, P_{nb,t}, Q_{nb,t}, v_{bt}, ~\mbox{s.t} \nonumber \\
    ~ & -r_{gt} \le \Delta P_{gt} \le r_{gt}, \forall g \in{\mathcal{G}} \label{eq:Fc.1}\\
    ~ & P_{bt}=P_{bt}^{0}-\sum\nolimits_{g \in C(b)}(P_{gt}+\Delta P_{gt})+\sum\nolimits_{k\in C( b)}q_{kt}\nonumber\\
    ~& ~~~ ~~~ -\sum\nolimits_{e\in C(b)}(p_{et}^{d}-p_{et}^{c}),\forall b\in{\mathcal{B}}, \label{eq:Fc.2}\\
    ~ & Q_{bt}=Q_{bt}^{0}-\sum\nolimits_{g \in C(b)}Q_{gt},\forall b\in{\mathcal{B}}, \label{eq:Fc.3}\\
    ~ & P_{nb,t}-P_{bt}-\sum\nolimits_{(b,m)\in{\mathcal{L}}} P_{bm,t}=0,\forall (n,b)\in{\mathcal{L}},\label{eq:Fc.4}\\
     ~ &Q_{nb,t}-Q_{bt}-\sum\nolimits_{(b,m)\in{\mathcal{L}}} Q_{bm,t}=0,\forall (n,b)\in{\mathcal{L}}, \label{eq:Fc.5}\\   
     ~&-\overline{P}_{nb}\le P_{nb,t} \le \overline{P}_{nb},\forall (n,b)\in{\mathcal{L}} ,\label{eq:Fc.6}\\ 
     ~&-\overline{Q}_{nb}\le Q_{nb,t} \le \overline{Q}_{nb},\forall (n,b)\in{\mathcal{L}},\label{eq:Fc.7}\\  ~&v_{nt}-(r_{nb}P_{nb,t}+x_{nb}Q_{nb,t})=v_{bt},\forall (n,b)\in{\mathcal{L}}, \label{eq:Fc.8}\\
     ~&\underline{v}_{b} \le {v}_{bt} \le  \overline{v}_{b},\forall b \in {\mathcal{B}}, ~{v}_{1t}=V_{0},\label{eq:Fc.9}\\
     ~& \underline{p}_{et}^{DA,c} \le p_{et}^{c} \le \overline{p}_{et}^{DA,c},\forall e \in \mathcal{E}, \label{eq:Fc.10}\\
     ~& \underline{p}_{et}^{DA,d} \le p_{et}^{d} \le \overline{p}_{et}^{DA,d},\forall e \in \mathcal{E} \} \label{eq:Fc.11}
\end{align}
\esq
Constraint \eqref{eq:Fc.1} restricts the adjusted output of gas-fired units to vary within the range of operational reserves. Constraints \eqref{eq:Fc.2}-\eqref{eq:Fc.9} constitute the linearized DistFlow model \cite{bai2017distribution}. Constraints \eqref{eq:Fc.10}-\eqref{eq:Fc.11} allow the real-time charging/discharging power of energy storage units to vary within the operational bounds decided in the day-ahead stage (will be introduced later). Denote the objective function and feasible set of problem \eqref{eq:sharing-oper} as $\Gamma_t^o(\lambda_t)$ and $\mathcal{X}_t^o(d_t,b_t)$, respectively.

\textbf{Step 4:} If $\|\lambda_t^{n}-\lambda_t^{n+1}\| \leq \epsilon$, let $\lambda_t^*=\lambda_t^{n+1}$ and go to \textbf{Step 5}; otherwise, let $n=n+1$ and go to \textbf{Step 2}.

\textbf{Step 5:} Each smart meter $k \in \mathcal{K}$ determines the optimal demand $d_{kt}^*$, and sharing quantity $q_{kt}^*$ based on $\lambda_{kt}^*$, and sends them back to each customer to execute. 

\emph{Remark on the objective function \eqref{eq:sharing-oper.1}:} The objective function \eqref{eq:sharing-oper.1} is an essential design of our proposed energy sharing mechanism. As we prove in Proposition \ref{Thm:prop-central} later, by using this objective function, we can guarantee the existence of a unique, social optimal energy sharing equilibrium. Since generalized Nash games are difficult to analyze, obtaining such a conclusion is not trivial but technically challenging. Moreover, we can provide an economic intuition of the designed objective function \eqref{eq:sharing-oper.1} in the case without energy storage, i.e. minimizing the objective function \eqref{eq:sharing-oper.1} is equivalent to minimizing the variance of all prosumers' energy sharing prices. This ensures a fair market outcome for all prosumers.

The energy sharing market involves complicated interactions among customers. The energy sharing market equilibrium is defined as follows, which no one can become better-off by unilaterally deviating from.

\begin{definition} (Energy Sharing Market Equilibrium) A profile $(d_t^*,b_t^*,q_t^*,\lambda_t^*)$ \footnote{$d_t^*$ refers to the collection of $d_{kt}^*,\forall k$. Similar for $b_t^*$, $q_t^*$, and $\lambda_t^*$.} is an \emph{energy sharing market equilibrium} if
\begin{align}
    \forall k \in \mathcal{K}:~ (d_{kt}^*,b_{kt}^*,q_{kt}^*) = ~ &  \mbox{argmin}_{d_{kt},b_{kt},q_{kt}} ~ \Gamma_k^c(d_{kt},q_{kt},\lambda_{kt}^*) \nonumber\\
    ~ & \mbox{s.t.} (d_{kt},b_{kt},q_{kt}) \in \mathcal{X}_{kt}^c(\lambda_{kt}^*)
\end{align}
and
$
    \lambda_t^* = \mbox{argmin}_{\lambda_t} \; \{ \Gamma_t^o(\lambda_t) ~|~ \lambda \in \mathcal{X}_t^o (d_t^*,b_t^*)\}
$.
\label{def1}
\end{definition}

According to the above definition, both the objective function and the action set of one customer depend on the strategies of other customers. Therefore, they play a generalized Nash game \cite{harker1991generalized}.
Generally, there is no guarantee for the existence and uniqueness of the equilibrium of a generalized Nash game. Fortunately, for the proposed energy sharing game, we are able to prove the existence of a partially unique market equilibrium, which is given in Section \ref{sec-3}. 

\subsection{Endogenous Uncertainty Set}

The real-time RG outputs $w_{kt},\forall k \in \mathcal{K}, \forall t \in \mathcal{T}$ vary within an uncertainty set $\mathcal{W}(u)$, depending on the RG connection/disconnection strategy $u$ in the day-ahead stage.
\bsq
\label{eq:DDU-set}
\begin{align}
    \mathcal{W}(u) =~& \{ w_{kt}, \forall k \in \mathcal{K}, \forall t \in \mathcal{T} ~| \nonumber \\
    &~\underline{W}_{kt}u_{kt} \le w_{kt} \le \overline{W}_{kt}u_{kt},\forall k \in \mathcal{J}, \forall t \in \mathcal{T} \label{eq:DDU-set.1}\\
    ~& \sum \nolimits_{k \in \mathcal{J}} \frac{|w_{kt}-w_{kt}^e|}{w_{kt}^h} \le B^S,\forall t \in \mathcal{T} \label{eq:DDU-set.2}\\
    ~ & \sum \nolimits_{t \in \mathcal{T}} \frac{|w_{kt}-w_{kt}^e|}{w_{kt}^h} \le B^T, \forall k \in \mathcal{J} \label{eq:DDU-set.3}\\
    ~ & w_{kt}^e = 0.5(\underline{W}_{kt}u_{kt}+\overline{W}_{kt}u_{kt}),\forall k \in \mathcal{J},\forall t \in \mathcal{T} \label{eq:DDU-set.4}\\
    ~ & w_{kt}^h = 0.5(\overline{W}_{kt}-\underline{W}_{kt}),\forall k \in \mathcal{J},\forall t \in \mathcal{T} \} \label{eq:DDU-set.5}
\end{align}
\esq
If RG $k \in \mathcal{J}$ is disconnected in period $t \in \mathcal{T}$, we have $u_{kt}=0$, and thus $w_{kt}=0$ according to \eqref{eq:DDU-set.1}. If RG $k \in \mathcal{J}$ is connected, its real-time power varies within a confident interval $[\underline{W}_{kt},\overline{W}_{kt}]$. 
Moreover, uncertainty budgets $B^S$ and $B^T$ are introduced to avoid over conservativeness. A disconnected RG does not contribute to the budget constraint \eqref{eq:DDU-set.2}-\eqref{eq:DDU-set.3} as both $w_{kt}$ and $w_{kt}^e$ are zero. 
The above uncertainty set depends on the day-ahead stage strategy $u$, and is an \emph{endogenous uncertainty set}. 

\begin{figure}[t]
\centering
\includegraphics[width=0.95\columnwidth]{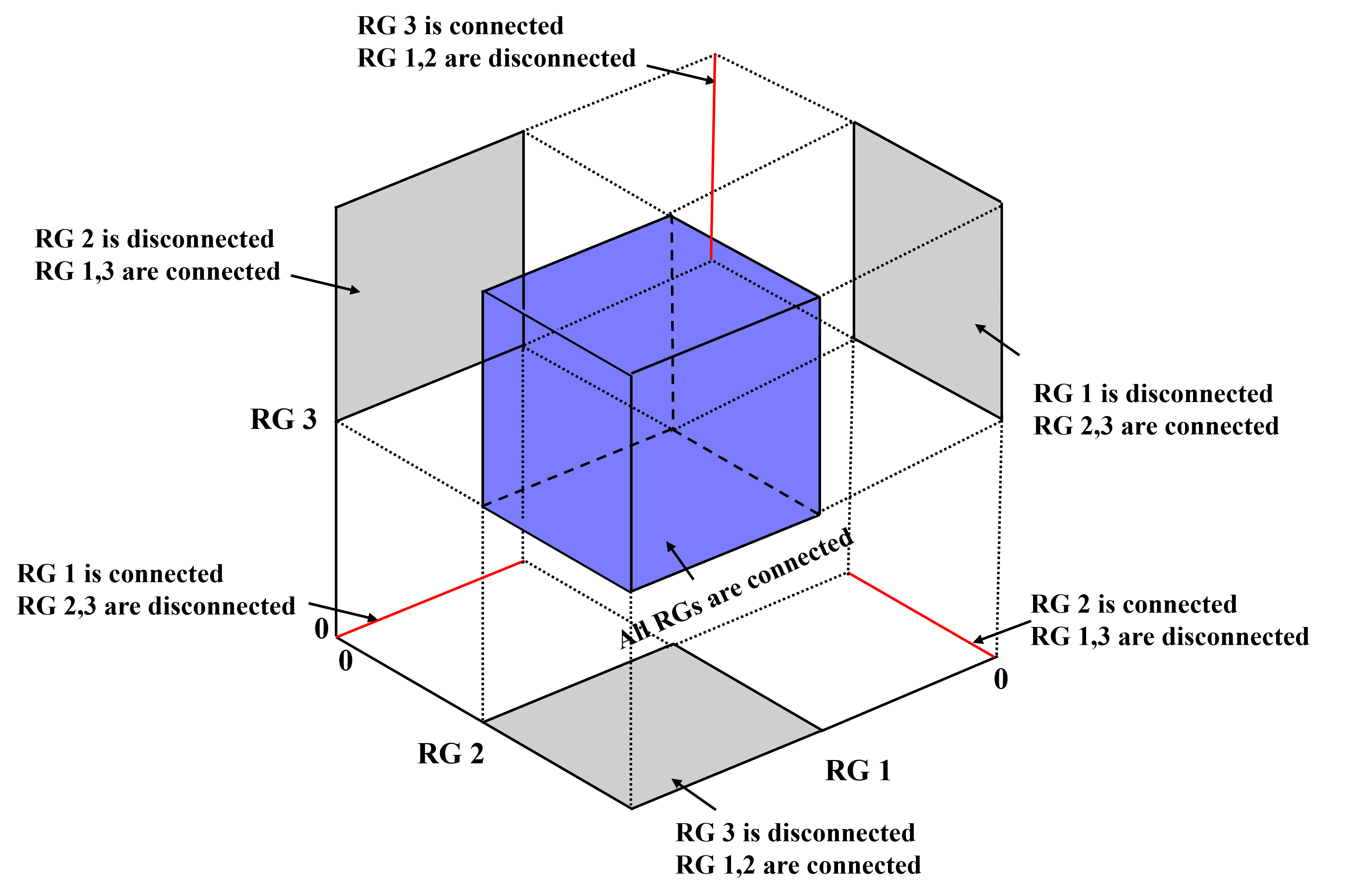}
\caption{Uncertainty sets under different disconnection strategies.}
\label{fig:US-UC}
\end{figure}

A simple example with three RGs is given in Fig. \ref{fig:US-UC} to illustrate how the RG connection/disconnection strategy influences the uncertainty set. If all RGs are connected, the uncertainty set will be the blue cube; if one of them is disconnected, the uncertainty set turns into a grey rectangle; if only one RG is connected, the uncertainty set is a red segment. Therefore, the dimension and shape of the uncertainty set change when the connection/disconnection strategy alters.

\emph{Remark on the uncertainty set:} The conservativeness of the uncertainty set can be adjusted by changing the budget parameters $B^S$ and $B^T$. If $B^S$ is set to be the number of prosumers $J$ and $B^T$ is set to be $T$, then the uncertainty set contains all the possible scenarios including those rare events. However, in practice, it is too conservative and costly to consider all the rare events with extremely small probabilities. By contrast, rare events are addressed by load shedding, generation re-dispatch and tripping, dynamic braking, and control system separation (or islanding) \cite{huang2019adaptive}.

\subsection{Two-Stage Robust Microgrid Dispatch}
\label{sec-2-2}
In the day-ahead stage, the microgrid operator decides on the connection/disconnection $u_{kt} ,\forall t \in \mathcal{T}$ of all RGs $k \in \mathcal{J}$, and the expected output and reserve of all dispatchable gas-fired units $P_{gt},Q_{gt},\forall g\in \mathcal{G},\forall t\in \mathcal{T}$, trying to minimize the total cost, including the penalty cost for underutilized renewable energy, the operation and reserve costs of gas-fired units, and the worst-case total disutility of all customers. Additionally, to facilitate real-time online operation of energy storage units, the bounds of their charging/discharging power and SOC levels are also optimized in the day-ahead stage. The operator faces a tradeoff between security and economy: If too many RGs are disconnected, more electricity needs to be generated by gas-fired units to meet the demand, leading to a higher cost. Otherwise, if few RGs are disconnected, the uncertainty increases and the system may collapse when the variation exceeds the demand adjustment capability.

With the real-time energy sharing game \eqref{eq:sharing-pro}-\eqref{eq:sharing-oper} and the endogenous uncertainty set \eqref{eq:DDU-set}, the microgrid operator solves the following robust optimization problem:
\bsq
\label{eq:robust-upper}
\begin{align}
    & \begin{aligned}
    \mathop{\min}_{x,u} ~ & \sum\nolimits_{k \in \mathcal{J}} \sum\nolimits_{t \in \mathcal{T}} \alpha_{kt}(1-u_{kt}) \\
    & + \sum\nolimits_{g\in\mathcal{G}} \sum\nolimits_{t\in\mathcal{T}} (c_g P_{gt}+s_g r_{gt}) + \mathop{\max}_{w \in \mathcal{W}(u)} F(x,w)
    \end{aligned}  \label{eq:robust-upper.1}\\
    & \mbox{s.t.} ~  \forall w \in \mathcal{W}(u),~ \exists ~\text{an equilibrium of game \eqref{eq:sharing-pro}-\eqref{eq:sharing-oper}} \\
    & x = (h_{gt}, P_{gt}, Q_{gt}, \mu_{et}^c, \mu_{et}^d, \underline{E}_{et}^{DA}, \overline{E}_{et}^{DA}, \underline{p}_{et}^{DA}, \overline{p}_{et}^{DA}, \forall g, \forall e, \forall t) \label{eq:robust-upper.x}\\
      & u_{kt} \in \{0,1\},\forall k \in \mathcal{J},\forall t \in \mathcal{T} \label{eq:robust-upper.2}\\
      & h_{gt} \in \{0,1\},\forall g \in \mathcal{G}, \forall t \in \mathcal{T} \label{eq:robust-upper.3}\\
      & h_{gt} \underline{P}_{g} \le P_{gt}+r_{gt},P_{gt}-r_{gt} \le h_{gt} \overline{P}_{g},\forall g \in \mathcal{G}, \forall t \in \mathcal{T} \label{eq:robust-upper.4}\\
       & h_{gt} \underline{Q}_{g} \le Q_{gt} \le h_{gt} \overline{Q}_{g},\forall g \in \mathcal{G}, \forall t \in \mathcal{T} \label{eq:robust-upper.5}\\
       & \mu_{et}^{c},\mu_{et}^{d} \in \{0,1\},\forall e \in \mathcal{E}, \forall t \in \mathcal{T} \label{eq:robust-upper.6}\\
       & \underline{E}_{et}^{DA}=\underline{E}_{e,t-1}^{DA}+(\underline{p}_{et}^{DA,c}\eta_{e}^{c}-\overline{p}_{et}^{DA,d}/\eta_{e}^{d})\Delta t,\forall e \in \mathcal{E}, \forall t \in \mathcal{T} \label{eq:robust-upper.7}\\
       & \overline{E}_{et}^{DA}=\overline{E}_{e,t-1}^{DA}+(\overline{p}_{et}^{DA,c}\eta_{e}^{c}-\underline{p}_{et}^{DA,d}/\eta_{e}^{d})\Delta t,\forall e \in \mathcal{E}, \forall t \in \mathcal{T} \label{eq:robust-upper.8}\\
      &\mu_{et}^{c}+\mu_{et}^{d}\le 1,\forall e \in \mathcal{E}, \forall t \in \mathcal{T} \label{eq:robust-upper.9}\\
       & \mu_{et}^{c}\underline{p}_{e}^{c} \le \underline{p}_{et}^{DA,c} \le \overline{p}_{et}^{DA,c} \le \mu_{et}^{c}\overline{p}_{e}^{c},\forall e \in \mathcal{E}, \forall t \in \mathcal{T} \label{eq:robust-upper.10}\\
       & \mu_{et}^{d}\underline{p}_{e}^{d} \le \underline{p}_{et}^{DA,d} \le \overline{p}_{et}^{DA,d} \le \mu_{et}^{d}\overline{p}_{e}^{d},\forall e \in \mathcal{E}, \forall t \in \mathcal{T} \label{eq:robust-upper.11}\\
       & \underline{E}_{e} \le \underline{E}_{et}^{DA} \le \overline{E}_{et}^{DA} \le \overline{E}_{e},\forall e \in \mathcal{E}, \forall t \in \mathcal{T} \label{eq:robust-upper.12}\\
        & -\Delta_{e} \le \underline{E}_{eT}^{DA}-E_{e0} \le \overline{E}_{eT}^{DA}-E_{e0} \le \Delta_{e},\forall e \in \mathcal{E} \label{eq:robust-upper.13}
\end{align}
\esq
where
\begin{align}
\label{eq:robust-lower}
    F(x,w)=\sum \nolimits_{k \in \mathcal{K}} \sum \nolimits_{t \in \mathcal{T}} U_k(d_{kt}^*(x,w))
\end{align}
and $(d^*(x,w),b^*(x,w),q^*(x,w),\lambda^*(x,w))$ is an energy sharing market equilibrium given $x$ and $w$. The objective function \eqref{eq:robust-upper.1} minimizes the total cost under the worst-case scenario, consisting of the underutilization cost of RGs, operation and reserve cost of gas-fired units, and the total cost of customers. It is worth noting that $F(x,w)$ in \eqref{eq:robust-lower} only consists of the total disutility. The total net payment in the energy sharing market is excluded since it is an internal trading cost not being concerned by the microgrid operator. We prove later in Proposition \ref{Thm:prop-payment} that the total net payment in the energy sharing market is non-negative, meaning that no subsidy from the microgrid operator is needed to support its operation. 

Constraint \eqref{eq:robust-upper.x} denotes all the decision variables in the day-ahead stage. RG connection/disconnection strategy related constraint is \eqref{eq:robust-upper.2}. \eqref{eq:robust-upper.3}-\eqref{eq:robust-upper.5} determine the on/off status, and limit the range of active/reactive power and operational reserves of controllable gas-fired units. Constraint \eqref{eq:robust-upper.6} determines the charged/discharged status of ESs. \eqref{eq:robust-upper.7}-\eqref{eq:robust-upper.8} are the constraints for the charging/discharging power and SOC level of energy storage units, to ensure the non-anticipativity and all-scenario feasibility in the real-time stage. The effectiveness of constraints \eqref{eq:robust-upper.9}-\eqref{eq:robust-upper.13} is proven in Proposition \ref{Thm:prop-SOC} below.

\begin{proposition} 
For each period $t \in \mathcal{T}$, suppose the charging power $p_{et}^c$ and discharging power $p_{et}^d$ satisfy \eqref{eq:Fc.10} and \eqref{eq:Fc.11}, respectively, where the parameters ($\underline{p}_{et}^{DA,c}$, $\overline{p}_{et}^{DA,c}$, $\underline{p}_{et}^{DA,d}$, $\overline{p}_{et}^{DA,d}$) are determined by \eqref{eq:robust-upper.9}-\eqref{eq:robust-upper.13}. Then, the SOCs of energy storage units always rest in the feasible physical ranges $[\underline{E}_e, \overline{E}_e],\forall e$, although these ranges are not explicitly considered in the real time problem \eqref{eq:sharing-oper}-\eqref{eq:Fc}.
\label{Thm:prop-SOC}
\end{proposition}
The proof of Proposition \ref{Thm:prop-SOC} can be found in Appendix \ref{apen-SOC}. It tells us that by setting non-anticipative operational ranges [$\underline{p}_{et}^{DA,c}$, $\overline{p}_{et}^{DA,c}$] 
 and [$\underline{p}_{et}^{DA,d}$, $\overline{p}_{et}^{DA,d}$] for charging and discharging power in the day-ahead stage, the energy storage units can operate in an online manner in the real-time stage.

The robust microgrid dispatch problem \eqref{eq:robust-upper}-\eqref{eq:robust-lower} is difficult to solve due to two reasons: 1) the model is a robust optimization problem with endogenous uncertainty that is computationally intractable in general cases; 2) the lower-level problem is a generalized Nash game instead of an optimization problem as traditional, whose equilibrium is hard to derive and analyze.

\section{Solution Approach}
\label{sec-3}

In this section, we develop an effective solution algorithm for the RMD model \eqref{eq:robust-upper}-\eqref{eq:robust-lower}. Specifically, properties of the lower-level energy sharing game are proven, based on which we can turn the original RMD model into a robust optimization problem with a conventional ``min-max-min'' structure. Then, a projection-based C\&CG algorithm is developed to solve the equivalent robust model with endogenous uncertainty.

\subsection{Transformation of the Energy Sharing Game}
\label{sec-3-1}
The generalized Nash game \eqref{eq:robust-upper}-\eqref{eq:robust-lower} at the lower level of the RMD model exerts great challenges in solving the problem. For the proposed energy sharing game in particular, we can prove the existence of a partial unique equilibrium that can be computed by a centralized counterpart as Proposition \ref{Thm:prop-central} shows. Based on this, an equivalent robust model with the conventional ``min-max-min'' structure can be derived.

\begin{proposition}(Existence and Partial Uniqueness)
\label{Thm:prop-central}
An energy sharing market equilibrium of the game \eqref{eq:sharing-pro}-\eqref{eq:sharing-oper} exists if and only if problem \eqref{eq:central-lower} is feasible. Moreover, if $(d_t^*,b_t^*,q_t^*,\lambda_t^*)$ is such an equilibrium, then $d_t^*$ is the unique optimal solution of problem \eqref{eq:central-lower}. There exists $\hat \eta_{kt}$ dual optimal in problem \eqref{eq:central-lower} such that $\lambda_{kt}^*=\hat \eta_{kt},\forall k \in \mathcal{K}$. $q_{kt}^* = d_{kt}^* + d_{kt}^f, \forall k \in \mathcal{I}$ and $q_{kt}^* = d_{kt}^* + d_{kt}^f - w_{kt}, \forall k \in \mathcal{J}$, $b_{kt}^* = q_{kt}^* + a \lambda_{kt}^*, \forall k \in \mathcal{K}$.
\bsq
\label{eq:central-lower}
\begin{align}
    \mathop{\min}_{d_t,q_t} ~ & \sum \nolimits_{k \in \mathcal{K}} U_k(d_{kt}), \label{eq:central-lower.1}\\
   \mbox{s.t.} ~ &  \left\{\begin{aligned}
     & q_{kt}=d_{kt}+d_{kt}^{f},~\mbox{if}~ k \in \mathcal{I} \\
     & w_{kt}+q_{kt}=d_{kt}+d_{kt}^{f},~\mbox{if}~ k \in \mathcal{J}
    \end{aligned}
    \right. : \eta_{kt}, \label{eq:central-lower.2}\\
    ~ & \underline{D}_{kt} \le d_{kt} \le \overline{D}_{kt},~ \forall k \in \mathcal{K},\label{eq:central-lower.3} \\
    ~ & (q_{kt},\forall k) \in \mathcal{F}_{c}. \label{eq:central-lower.4}
\end{align}
\esq
\end{proposition}

The proof of Proposition \ref{Thm:prop-central} can be found in Appendix \ref{apen-central}. When the robust optimization problem is solved, problem \eqref{eq:central-lower} always has a feasible point. Therefore, there always exists an energy sharing market equilibrium and the quantity $d_t^*$ is unique. This ensures the practicability of the proposed mechanism. Moreover, Proposition \ref{Thm:prop-central} offers a centralized counterpart to compute the equilibrium. Noticing that the objective \eqref{eq:central-lower.1} is the same as $F(x,w)$ in \eqref{eq:robust-lower}. Therefore, the original lower-level game \eqref{eq:sharing-pro}-\eqref{eq:sharing-oper} can be equivalently replaced by optimization \eqref{eq:central-lower}. By doing so, the original robust model turns into a ``min-max-min'' structure as follows.
\bsq
\label{eq:robust-eq}
\begin{align}
    \mathop{\min}_{x,u} ~ & \sum\nolimits_{k \in \mathcal{J}} \sum\nolimits_{t \in \mathcal{T}} \alpha_{kt}(1-u_{kt}) + \sum\nolimits_{g\in\mathcal{G}} \sum\nolimits_{t\in\mathcal{T}} (c_g P_{gt}+s_g r_{gt}) \nonumber\\
    & + \mathop{\max}_{w \in \mathcal{W}(u)} \min_{d \in \mathcal{X}^d(x,w)} \sum\nolimits_{k \in \mathcal{K}} \sum\nolimits_{t \in \mathcal{T}} U_k(d_{kt}) \label{eq:robust-eq.1}\\
    & \mbox{s.t.} ~ \text{\eqref{eq:robust-upper.x}-\eqref{eq:robust-upper.13}}, \mathcal{X}^d(x,w) \neq \varnothing, \forall w \in \mathcal{W}(u)
\end{align}
\esq
where $\mathcal{X}^d(x,w) = \{d~|~ \exists ~q ~\mbox{satisfies}~\text{\eqref{eq:central-lower.2}-\eqref{eq:central-lower.4}}\}$.

In fact, the above model \eqref{eq:robust-eq} is the RMD model under a centralized scheme, where the microgrid operator monitors all the consumers and prosumers within it. Therefore, the economy and flexibility of the proposed RMD model with real-time energy sharing are the same as that of the RMD model under centralized operation. It means that the total cost under the worst-case scenario is the same; if a renewable power output scenario can be balanced by centralized dispatch in real-time, then given this renewable output, there exists a corresponding energy sharing game equilibrium. Apart from the existence of an equilibrium, another issue the microgrid operator cares about is whether it needs to offer financial support for the energy sharing market operation, which is discussed in the proposition below.

\begin{proposition}\label{Thm:prop-payment} Suppose $(d_t^*,b_t^*,q_t^*,\lambda_t^*)$ is an equilibrium of the energy sharing game, then the total net payment of all customers in the energy sharing market is non-negative, i.e.
\begin{align}
    \sum \nolimits_{k \in \mathcal{K}}  \lambda_{kt}^* q_{kt}^* \ge 0
\end{align}
\end{proposition}
The proof of Proposition \ref{Thm:prop-payment} is in Appendix \ref{apen-payment}. It indicates that the total net payment of all customers in the energy sharing market is zero or positive, and thus, no external subsidy from the microgrid operator is needed. The positive net payment can be used to support the operation of the energy sharing market, maintenance of the trading platform, etc. 

\subsection{Linearization of the Objective Function}
\label{sec-3-2}
\begin{figure}[t]
\centering
\includegraphics[width=0.7\columnwidth]{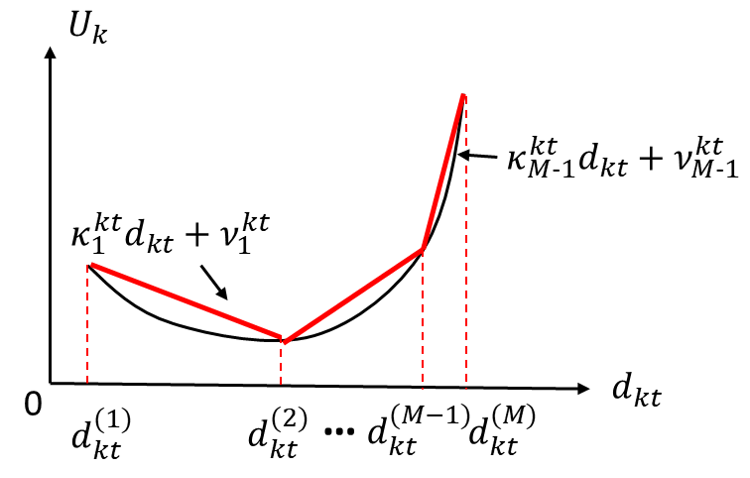}
\caption{Illustration of the objective function linearization.}
\label{fig:linearization}
\end{figure}

After the above transformation, the only nonlinearity of \eqref{eq:robust-eq} exists in the objective function. Since the disutility function $U_k(d_{kt})$ is strictly convex, it can be linearized using piecewise linearization technique. As demonstrated in Fig. \ref{fig:linearization}, we replace $U_k(d_{kt})$ with $\sigma_{kt}$, and add constraints
\bsq
\begin{align}
& \sigma_{kt} \ge \kappa_{m}^{kt} d_{kt}+ \nu_{m}^{kt},\forall m=1,...,M-1, \\
& d_{kt}^{(1)} \le d_{kt} \le d_{kt}^{(M)},
\end{align}
\esq
where $M$ is the number of sample points. Therefore, the original RMD model \eqref{eq:robust-upper}-\eqref{eq:robust-lower} turns into a ``min-max-min'' model with linear objective functions and constraints in all stages. The compact form is
\bsq
\label{eq:compact}
\begin{align}
    \min_{x,u} ~ & \gamma^\top x + \beta^\top u + \max_{w \in \mathcal{W}(u)} \min_{y \in \Pi(x,w)} \rho^\top y, \\
    \mbox{s.t.}~ & Ax \ge e, u \in \{0,1\}^{|\mathcal{J}| \times |\mathcal{T}|},
\end{align}
\esq
where
\begin{align}
    & \mathcal{W}(u)=\{w ~|~\exists w',~ \mbox{s.t.}~ R w + A w' \le g-Vu\}, \nonumber \\
    & \Pi(x,w)=\{y ~|~ H y \ge h-C x - T w\} .\nonumber
\end{align}

Different from the conventional robust optimization model where the uncertain factor $w$ varies within a fixed set $\mathcal{W}$, here the set $\mathcal{W}(u)$ depends on the day-ahead decision $u$. Typical methods to solve the conventional robust optimization with exogenous uncertainty, such as C\&CG, cannot be directly applied due to two reasons: 1) The previously selected scenario may fall outside of the uncertainty set when the day-ahead decision $u$ changes. For example, in Fig. \ref{fig:US-UC}, if RG 1 is connected in the current iteration and then disconnected in the next iteration with RGs 2-3 connected during both iterations, then the scenario selected in the current iteration is inside the blue cube but outside the uncertainty set in the next iteration (the grey rectangle). 2) In problem \eqref{eq:compact}, the day-ahead strategy $u$ only influences the uncertainty set $\mathcal{W}(u)$ but does not appear in the real-time feasible set $\Pi(x,w)$. Therefore, when applying C\&CG, adding the worst-case scenario into the master problem will result in an unchanged day-ahead strategy $u$ and the algorithm gets stuck. To tackle the two issues above, we develop a projection-based C\&CG algorithm below that projects the previously selected scenarios into the new uncertainty set with a different $u$.

\subsection{Projection-based C\&CG Algorithm}
\label{sec-3-3}
First, take the dual problem of the inner ``min'' problem and combine it with the middle ``max'', we have
\bsq
\label{eq:inner-eq}
\begin{align}
   \max_{w,z} ~ & (h- C x - T w)^\top z, \\
    \mbox{s.t.}~ & z \in Z=\{z~|~ H^\top z = \rho, z \ge 0\}, \\
    ~ & w \in \mathcal{W}(u).
\end{align}
\esq

Due to the bilinear term $w^\top z$ in the objective function, problem \eqref{eq:inner-eq} is a nonconvex problem. Fortunately, the worst-case scenario always appears at the vertices of the uncertainty set \cite{lorca2014adaptive}, enabling us to equivalently substitute the uncertainty set \eqref{eq:DDU-set} by the special uncertainty set with binary variables. The mixed integer programming (MIP) method \cite{jiang2011robust} is used to gain the global optimal solution, and the big-M method is adopted to linearize the bilinear products of binary variables and continuous variables.
The alternating direction method \cite{jiang2011robust} is another common approach to deal with the problem \eqref{eq:inner-eq}. As later shown in Section \ref{sec-4-1}, the MIP method gains global optimality, whereas the alternating direction method may lead to sub-optimality. The merit of the alternating multiplier method is explicit in large-scale problems with faster convergence speed and acceptable optimality gap.

   However, the problem \eqref{eq:inner-eq} may be infeasible for a given ($x$,$u$). To check feasibility, we construct the following feasibility-check (\textbf{FC}) problem, which is also solved using inner problem duality and the alternating direction method.
\bsq
\label{eq:feasibility-eq}
\begin{align}
   \textbf{FC:}~ \max_{w \in \mathcal{W}(u)} \min_{y \in \Pi(x,w),s} ~ & 1^\top s ,\\
    \mbox{s.t.}~ & H y+s \ge h-C x - T w, s \ge 0.
\end{align}
\esq
   
\begin{figure}[t]
\centering
\includegraphics[width=0.9\columnwidth]{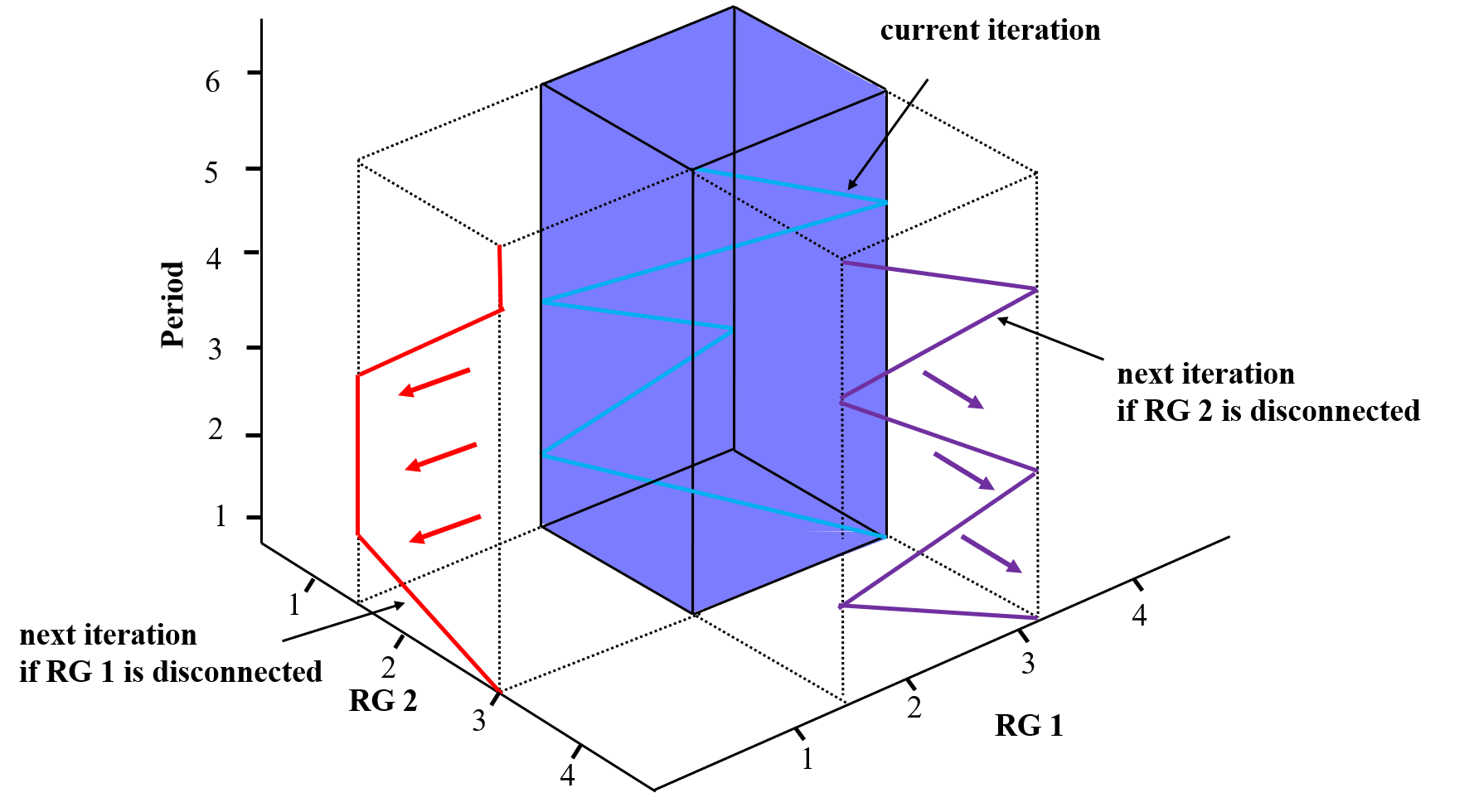}
\caption{Illustration of the projection operation.}
\label{fig:projection}
\end{figure}

Given the day-ahead stage decision $(x,u)$, the feasibility-check problem \eqref{eq:feasibility-eq} is first solved to check if the problem  \eqref{eq:inner-eq} is feasible. If not, a feasibility cut will be generated; else, we solve the problem \eqref{eq:inner-eq} to get an optimality cut. Conventionally, the worst-case scenarios are added to the master problem directly to generate feasibility/optimality cuts. However, due to the two difficulties caused by  endogenous uncertainty mentioned in Section \ref{sec-3-2}, adding the scenarios directly can lead to an over-conservative solution or infeasibility. To deal with this problem, we propose a modified master problem as follows with the previously selected scenarios projecting into the uncertainty set with a new $u$.
\bsq
\label{eq:master}
\begin{align}
    \textbf{Master:}~ \min_{x,u} ~ & \gamma^\top x + \beta^\top u + \tau, \\
    \mbox{s.t.}~ & Ax \ge e, u \in \{0,1\}^{|\mathcal{J}| \times |\mathcal{T}|}, \\
    ~ & H y^i \ge h- C x - T (u \circ w^i) ,\forall i \in \mathcal{S}_{w},\\
    ~ & \tau \ge \rho^\top y^i,\forall i \in \mathcal{S}_{w},
\end{align}
\esq
where $\mathcal{S}_{w}$ is the set of selected scenarios and $i$ is its index. For any element $i \in \mathcal{S}_w$, the projection operation $u \circ w^i$ turns $w^i$ into $u_{kt}w_{kt}^i,\forall k \in \mathcal{K},\forall t \in \mathcal{T}$, where $u_{kt}$ is the day-ahead decision variable to indicate whether RG $k$ is connected in period $t$ or not. Specifically, if $u_{kt}=1$, then after the operation $u \circ w^i$, we have $u_{kt}w_{kt}^i=w_{kt}^i \in [\underline{W}_{kt},\overline{W}_{kt}]$. If we set $u_{kt}=0$, then after the operation $u \circ w^i$, we have $u_{kt}w_{kt}^i=0$. Therefore, the projected scenario is always inside the uncertainty set. Fig. \ref{fig:projection} provides an intuitive illustration. The worst-case scenario selected previously is projected according to the current value of the day-ahead decision variable $u$. The modified master problem \eqref{eq:master} is still a mixed integer linear program (MILP) and is tractable. The projection-based C\&CG algorithm is summarized in Algorithm 1.

\begin{algorithm}[t]
        \caption{Projection-based C\&CG}
        \LinesNumbered 
        \KwIn{Initial set $\mathcal{S}_{w}=\{{w^e}\}$; tolerance $\epsilon > 0$; $j=1$; assign a large number to $UB_{0}$.}
        \KwOut{The RMD strategy $(x^*,u^*)$.}
        Solve the master problem \eqref{eq:master}; record the optimal solution $(x^j,u^j)$ and the optimal value $LB_{j}$\; 
        Solve the feasibility-check problem \eqref{eq:feasibility-eq} with $(x^j,u^j)$; if the optimal objective $>$ 0, let  $w^j$ be the worst scenario, let $UB_{j}\leftarrow UB_{j-1}$, and go to step 4; else, go to step 3\;
        With $(x^j,u^j)$, solve the problem \eqref{eq:inner-eq}; record the optimal solution $(w^j,z^j)$ and the optimal value $ub$; let $UB_{j}\leftarrow ub+\gamma^\top x^j+\beta^\top u^j$\;
        \eIf {$|UB_{j}-LB_{j}|\le \epsilon$}{$x^*\leftarrow x^j$; $u^*\leftarrow u^j$; terminate.}{ $\mathcal{S}_{w} \leftarrow \mathcal{S}_{w} \cup \{w^j\}$; $j\leftarrow j+1$; go to \textbf{Step 1}.}
\end{algorithm}

{\color{black}
\begin{proposition}\label{Thm:prop-converge} (Convergence of Algorithm 1) Assume the proposed two-stage RO with DDU model \eqref{eq:compact} has an optimal solution. Let $n := |\cup_{u \in \{0, 1\}^{|\mathcal{J}| \times |\mathcal{T}|}} V(\mathcal{W}(u))|$ be the total number of possible vertices of $\mathcal{W}(u)$ as $u$ varies in $\{0, 1\}^{|\mathcal{J}| \times |\mathcal{T}|}$. Then the projection-based C\&CG algorithm is guaranteed to converge to the optimal solution with an error of at most $\varepsilon$ within $\mathcal{O}(n)$ iterations.\end{proposition}

The proof of Proposition \ref{Thm:prop-converge} is in Appendix D. The main idea is as follows: The worst-case scenario always appears at a vertex of the uncertainty set in each iteration \cite{lorca2014adaptive}. Although the uncertainty set may change with the day-ahead decision variable $u$, $n$ is finite since $u$ can only take a finite number of possible values, and the number of vertices of uncertainty sets is also finite. Since a vertex is added into $\mathcal{S}_w$ in each iteration, Algorithm 1 is guaranteed to stop within $\mathcal{O}(n)$ iterations.
}

An alternative to Algorithm 1 is to turn \eqref{eq:robust-upper}-\eqref{eq:robust-lower} into an equivalent RO with exogenous uncertainty and solve it using the traditional C\&CG algorithm \cite{zhao2014variable} (we call it ``DIU-based method''). This can be achieved by  letting $w_{kt}^{'}=w_{kt}/u_{kt}$. However, as we show later in Section \ref{sec-4-1}, our proposed method is more efficient since it requires fewer iterations and auxiliary variables, and shorter computation time. Moreover, though only endogenous uncertainty is considered in this paper, our model and approach can be extended to incorporate both exogenous and endogenous uncertainties. Specifically, for exogenous uncertainty, we add the worst-case scenarios directly to the master problem; for endogenous uncertainty, we project the scenario onto the new uncertainty set using the approach in this paper. 

\emph{Remark on the difference between the projection-based C\&CG algorithm and ADMM algorithm:} The proposed projection-based C\&CG algorithm is completely different from the ADMM algorithm. For the mathematical model, the projection-based C\&CG algorithm is for solving two-stage robust optimization models, whereas the ADMM algorithm  is for solving a convex centralized optimization model in a distributed manner. For the solution procedure, the C\&CG algorithm divides the “min-max-min” problem into a master problem (the outer “min”) and a subproblem (the inner “max-min”), which will be solved iteratively. The master problem is solved to update the first-stage decision. Subsequently, the subproblem is solved under a fixed first-stage decision to identify the worst-case scenario, which is returned to the master problem. Then the master problem is solved again. This happens until convergence. By contrast, the ADMM algorithm starts by constructing the augmented Lagrangian function and iteratively updates the decision variables and dual variables until convergence.

\section{Case Studies}
\label{sec-4}
  
Numerical experiments are presented in this section to show the effectiveness, advantages, and scalability of the proposed model and algorithm.

 \subsection{Benchmark Case}
 \label{sec-4-1}

\begin{figure}[ht]
\centering
\includegraphics[width=1.02\columnwidth]{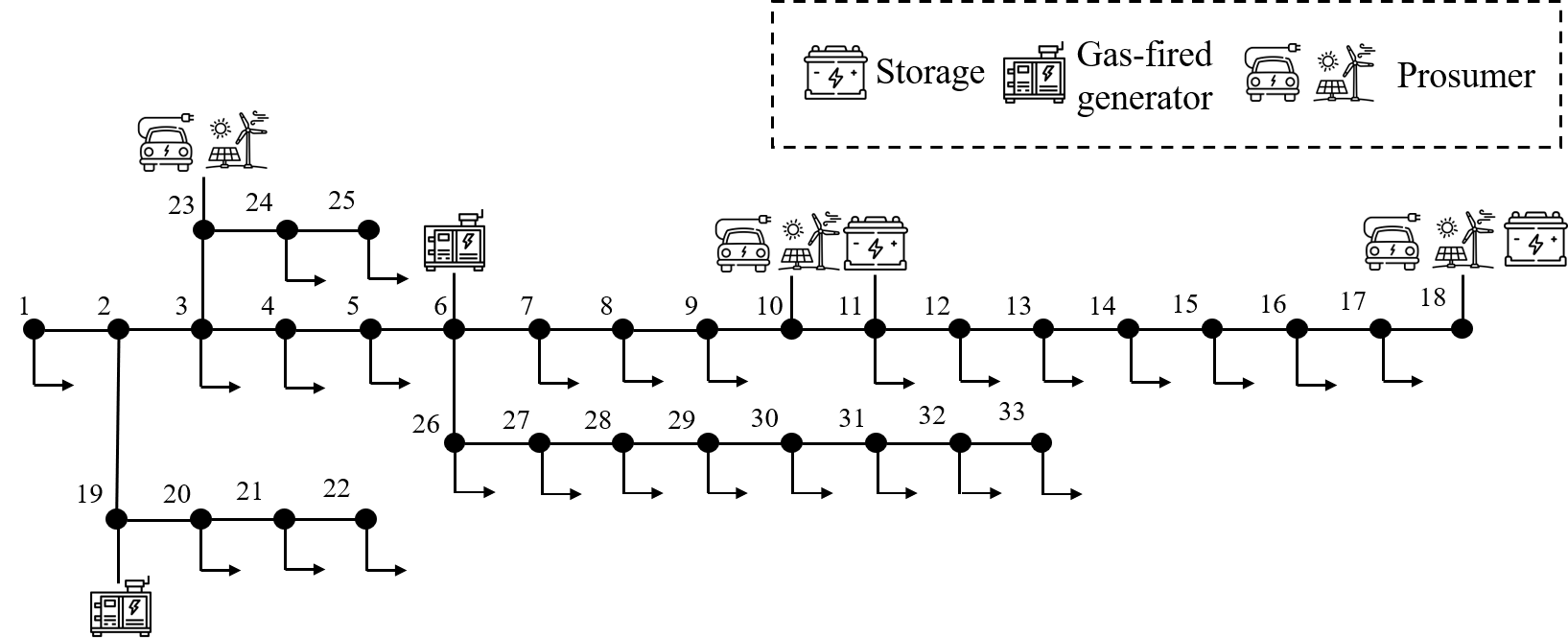}
\caption{The modified 33-bus microgrid system.}
\label{fig:33bus system}
\end{figure}

We first use the modified standalone 33-bus microgrid with three prosumers, two gas-fired generators, and two energy storage units for demonstration, as depicted in Fig. \ref{fig:33bus system}. Here, quadratic disutility functions $U_{kt}(d_{kt}):=\alpha_{k}^{1}(d_{kt}^2)-\alpha_{k}^{2}d_{kt}+\alpha_{k}^{3},\forall k$ are used with the parameters given in Table \ref{tab:par-prosumers} and the market sensitivity factor $a$ is set as 0.01 MW/\$. The penalty cost of underutilized renewable power output is 400 \$/MWh. The uncertainty budgets are $B^{S}=2$ and $B^{T}=4$.

\begin{table}[ht]
\renewcommand{\arraystretch}{1.3}
\renewcommand{\tabcolsep}{0.5em}
\centering
\caption{Parameters of prosumers}
\label{tab:par-prosumers}
\begin{tabular}{ccccccc} 
\hline
Prosumer & Bus & $\alpha_{k}^{1}$  &  $\alpha_{k}^{2}$ & 
$\alpha_{k}^{3}$ & $d_{kt}^f$ & $[\underline{D}_{kt},\overline{D}_{kt}]$\\
 & & \$/MW$^2$   & \$ /MW & \$ & MW & MW \\
\hline
1 & 10 & 30 & 360 & 1200 & 0.1 & $[0.1,3]$ \\
2 & 18 & 50 & 500 & 2300 & 0.2 & $[0.2,4]$\\
3 & 23 & 60 & 600 & 3000 &  0.15 & $[0.3,5]$ \\
\hline
\end{tabular}
\end{table}

\subsubsection{Energy sharing mechanism}

\begin{figure}[ht]
\centering
\subfigure[Energy sharing prices]{
\includegraphics[width=4.2cm,height=3.8cm]{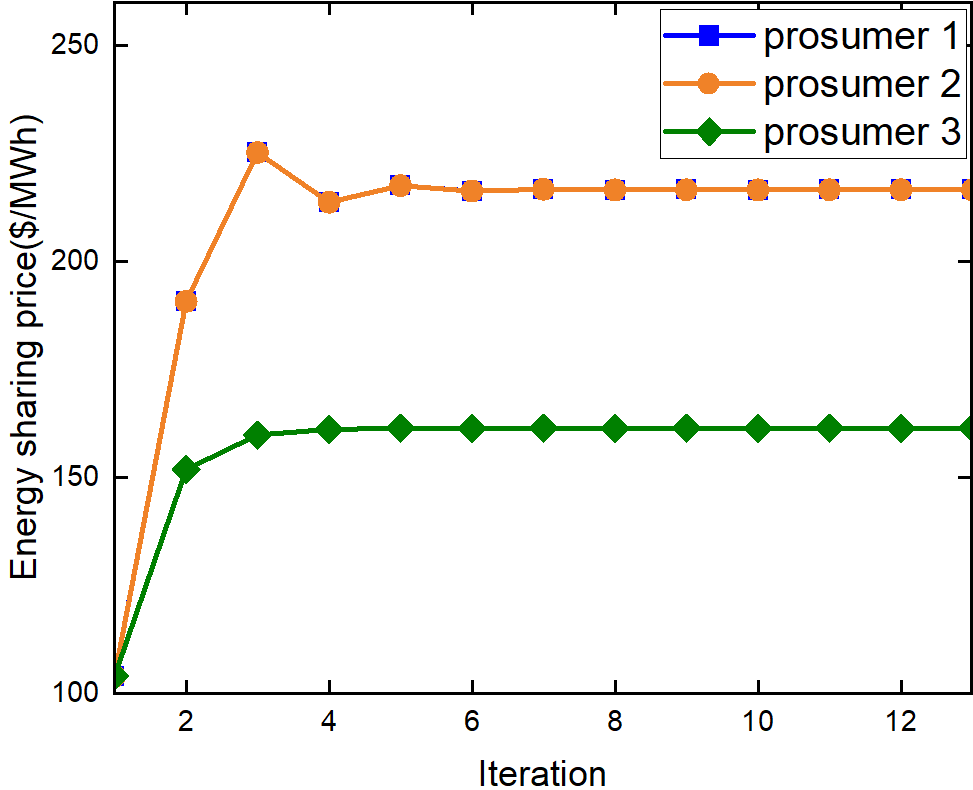}}\subfigure[Elastic demands]{
\includegraphics[width=4.2cm,height=3.8cm]{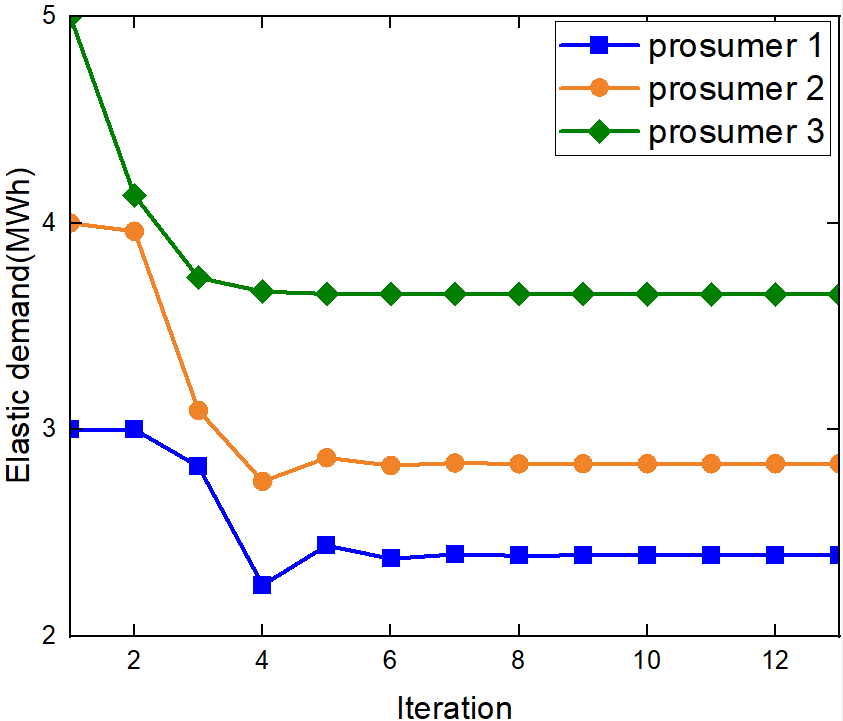}}
\caption{Energy sharing prices and elastic demands in each iteration.}
\label{fig:BP}
\end{figure}

\begin{figure}[ht]
\centering
\includegraphics[width=0.55\columnwidth]{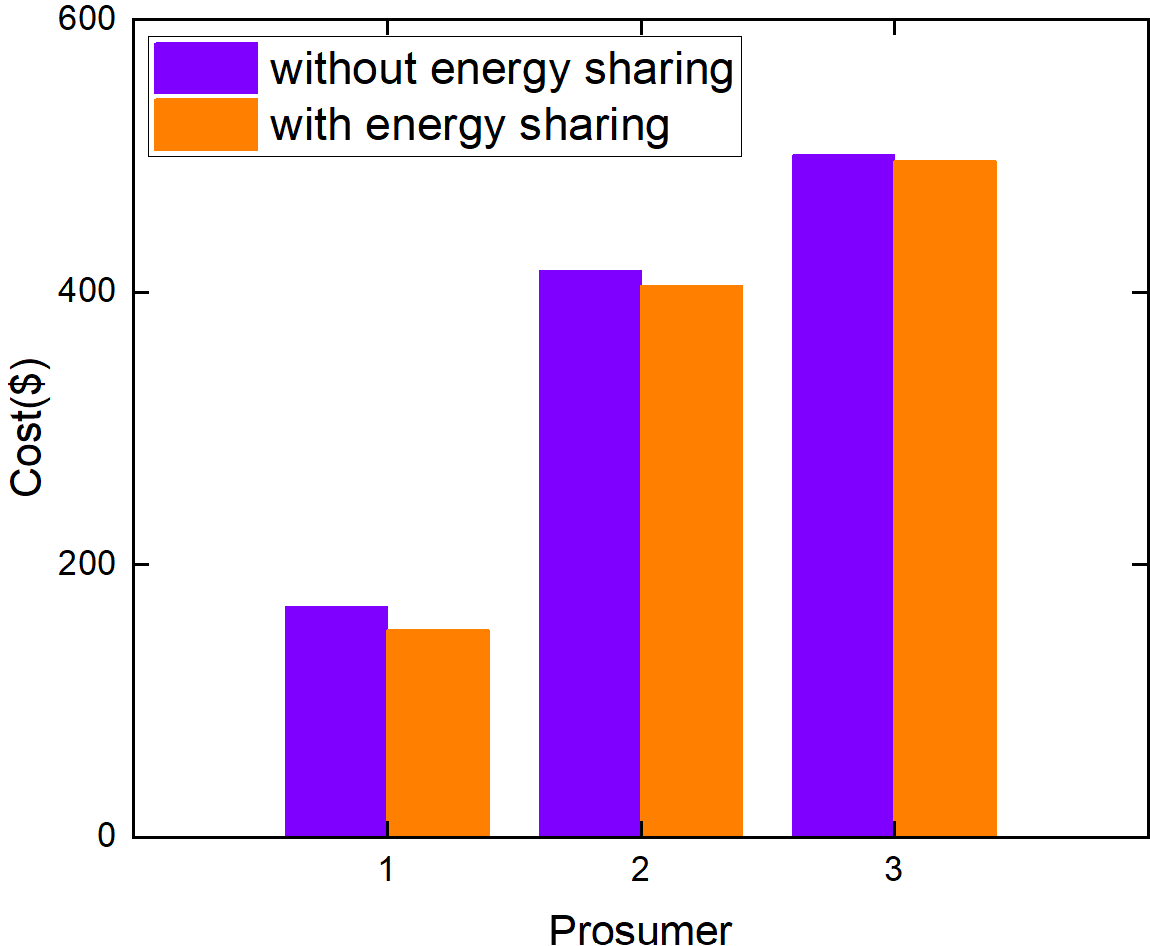}
\caption{Prosumers' cost comparison with and without energy sharing.}
\label{fig:cost}
\end{figure}

\begin{figure}[ht]
\centering
\subfigure[Energy sharing price]{
\includegraphics[width=4.2cm,height=3.8cm]{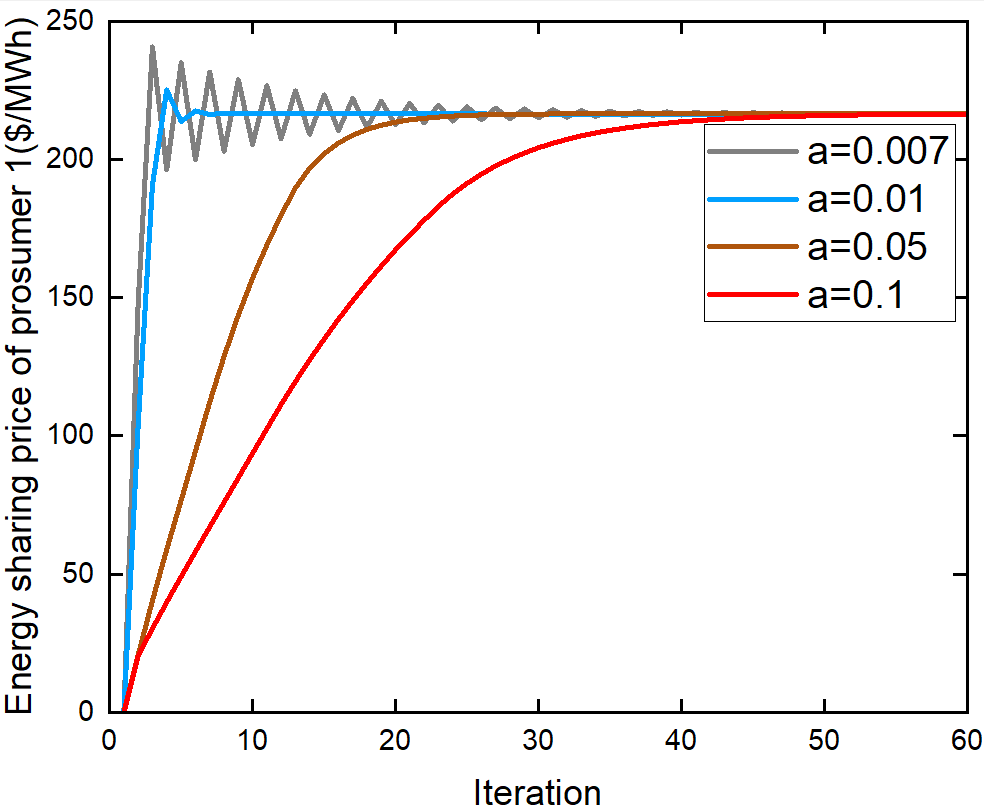}}\subfigure[Elastic demand]{
\includegraphics[width=4.2cm,height=3.8cm]{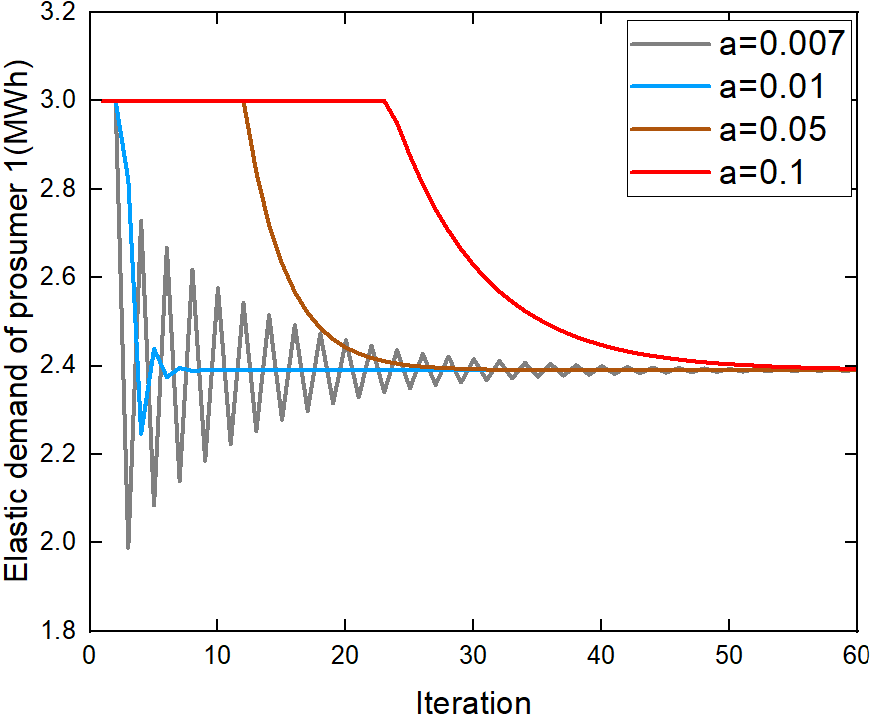}}
\caption{Impact of $a$ on prosumer 1's sharing price and elastic demand.}
\label{fig:a}
\end{figure}

\begin{figure}[ht]
\centering
\subfigure[Positive and negative energy sharing prices coexist]{
\includegraphics[width=4.2cm,height=3.8cm]{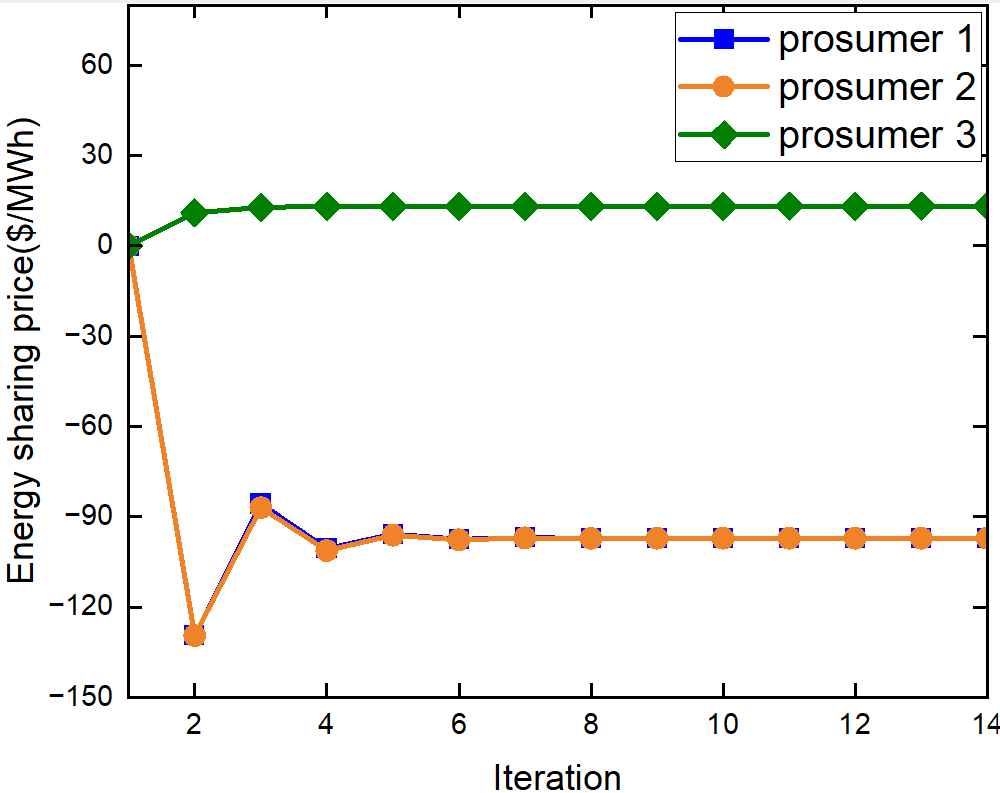}}\subfigure[Negative energy sharing prices]{
\includegraphics[width=4.2cm,height=3.8cm]{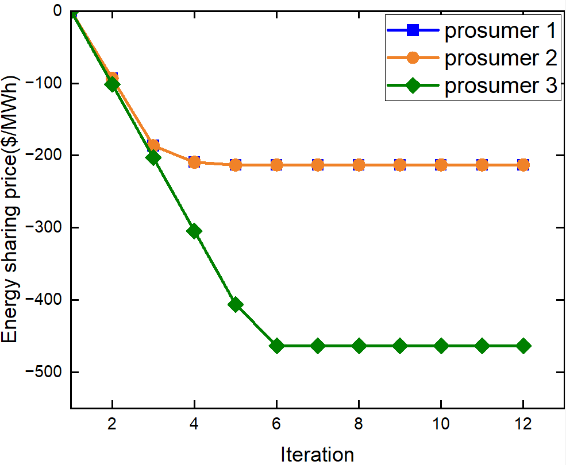}}
\caption{Cases with negative energy sharing prices.}
\label{fig:negative}
\end{figure}

First, we validate the effectiveness of the proposed energy sharing mechanism. Taking period 1 as an example, the changes of the prosumers' energy sharing prices and elastic demands are plotted in Fig. \ref{fig:BP}. We can observe that the energy sharing prices converge to 217, 217, 161 \$/MWh within 13 iterations, which equal the Lagrange multipliers of \eqref{eq:central-lower.2} at the optimum. The elastic demands converge to 2.39, 2.83, and 3.66 MW, respectively, which are the same as the optimal solution of \eqref{eq:central-lower}. The above validates Proposition \ref{Thm:prop-central}, which says that the energy sharing outcome is centralized optimal. Further, to show the advantages of the proposed energy sharing mechanism, we compare the costs of each prosumer between the cases with and without energy sharing. As shown in Fig. \ref{fig:cost}, the total cost of each prosumer decreases after sharing energy with other prosumers. Therefore, all prosumers are willing to participate in energy sharing. The proposed energy sharing mechanism reduces the total cost of the three prosumers by 3.00\%. In addition, we test the impact of market sensitivity $a$ on market equilibrium. Let $a = 0.007,0.01,0.05,0.1$, respectively, Fig.\ref{fig:a} shows the iteration process of prosumer 1's energy sharing price and elastic demand. We can observe that, although the convergence speeds are different, the algorithm converges to the identical market equilibrium under different $a$ as claimed in Proposition \ref{Thm:prop-central}. Specifically, a smaller $a$ leads to faster convergence speed, whereas too small $a$ may result in oscillation. Also, the sharing price could be negative in some cases, which is strongly related to the parameters of energy storage, DER, and disutility function. A negative sharing price arises when selling electricity results in a significant reduction in disutility. Consequently, the prosumer is willing to pay for selling electricity. We change $\alpha_{k}^{2}$ of the three prosumers to -30, -90, -600 \$/MW and 10, 0, 20 \$/MW, respectively. Correspondingly, Fig.\ref{fig:negative} shows the cases when the positive and negative energy sharing prices coexist and all prosumers' energy sharing prices are negative.

\subsubsection{Real-time energy storage operation}

\begin{figure}[ht]
\centering
\subfigure[Charging/discharging power bounds]{
\label{Fig.sub-1}
\includegraphics[width=4.2cm,height=3.8cm]{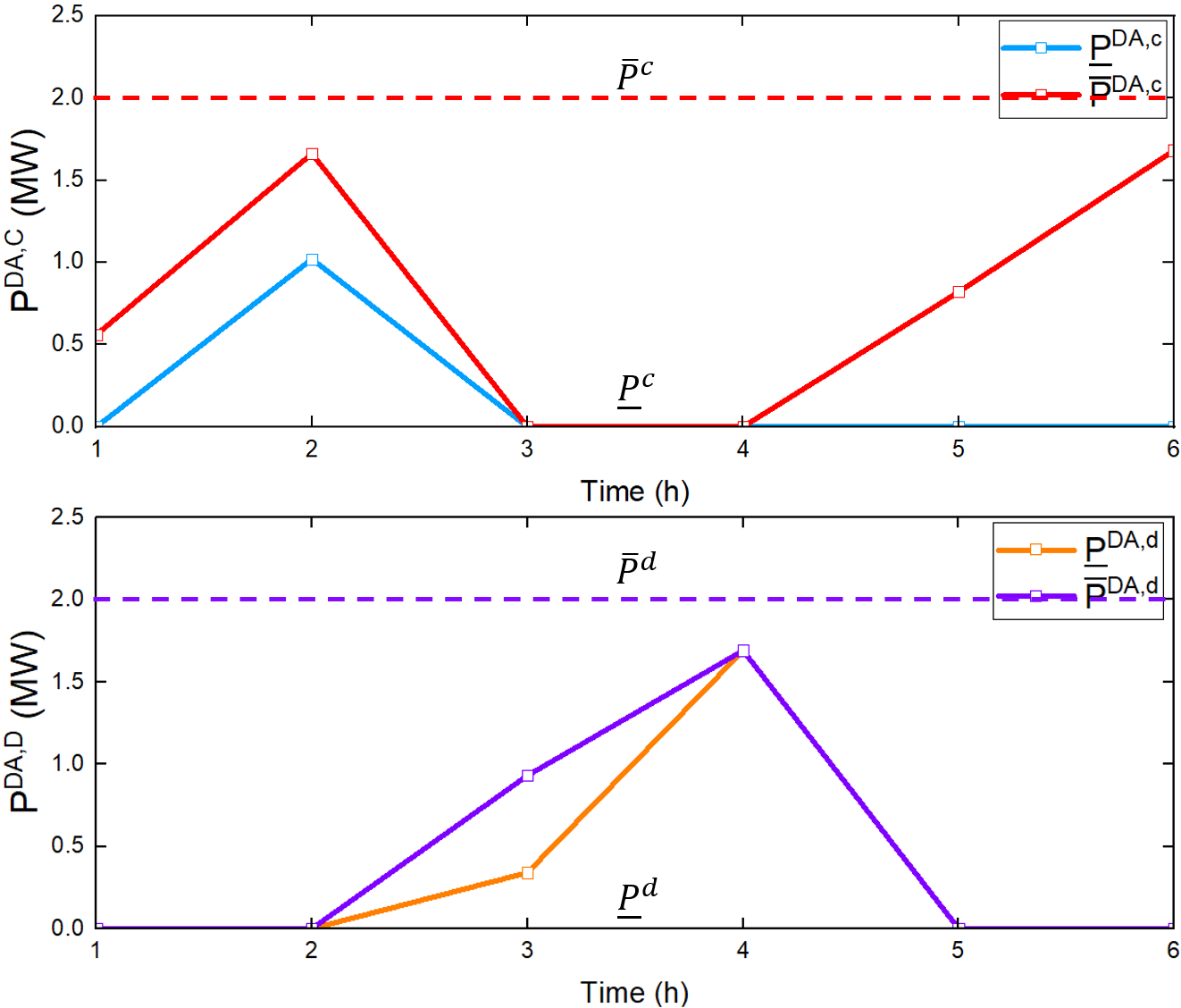}}\subfigure[SOC bounds and actual SOC, the expected renewable output]{
\label{Fig.sub-2}
\includegraphics[width=4.2cm,height=3.8cm]{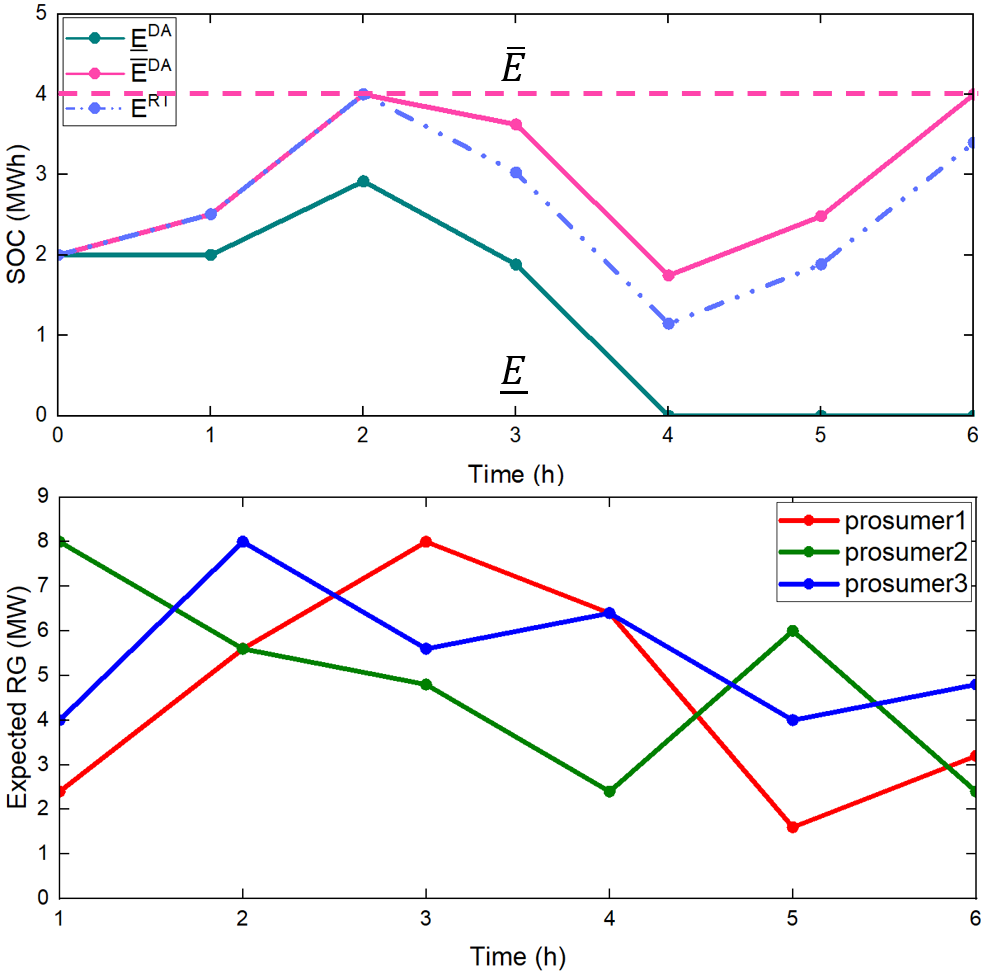}}
\caption{Operational bounds of energy storage units.}
\label{fig:es}
\end{figure}

We further investigate the impact of energy storage units on the dispatch strategy. 
Fig. \ref{Fig.sub-1} shows the charging/discharging power bounds of the energy storage unit at the node 11. It verifies that all the ES operational bounds $[\underline{P}_{et}^{DA,c}, \overline{P}_{et}^{DA,c}]$ and $[\underline{P}_{et}^{DA,d}, \overline{P}_{et}^{DA,d}]$ are subintervals of the technical intervals (the dash lines). Energy storage units store electricity in periods with high RG power outputs (e.g., periods 1 and 2) to avoid disconnection and discharge in periods when the RG power outputs are too low (e.g. periods 3 and 4). Fig. \ref{Fig.sub-2} demonstrates that when the charging and discharging power of energy storage satisfy \eqref{eq:Fc.10} and \eqref{eq:Fc.11}, the real-time SOC level $E^{RT}$ (the blue dot line) always falls in the pre-dispatched range $[\underline{E}^{DA},\overline{E}^{DA}]$ (the green and pink solid lines). Moreover, the strategic schedule of energy storage units helps to reduce the total underutilized RG output from 32.8 MWh to 21.6 MWh. Correspondingly, the total cost \eqref{eq:robust-upper.1} declines from 64,944\$ to 61,031\$.

\subsubsection{Projection-based C\&CG algorithm}
We first compare the computation performance between the alternating direction and the MIP method. We run the bilinear programming \eqref{eq:inner-eq} in 50 different scenarios by the two methods, respectively. The average computation time, optimal objective value, and optimality gap are shown in Table \ref{tab:par-bilinear}. The alternating multiplier method converges with an average time of 4.48s, while the MIP method needs an average time of 2.19s. The average optimality gap of the alternating multiplier method is 0.15\% compared to the exact value of bilinear programming attained by the MIP method, which is acceptable.

\begin{table}[ht]
\renewcommand{\arraystretch}{1.3}
\renewcommand{\tabcolsep}{0.5em}
\centering
\caption{Performance comparison between bi-linear algorithms}
\label{tab:par-bilinear}
\begin{tabular}{cccc} 
\hline 
  & Average   &  Average  & Average  \\
Method  & computation  &  optimal &  optimality  \\
  & time(s)  &  objective & gap  \\
\hline
Alternating direction & 4.48 &  29663.29 & 0.15\%\\
MIP & 2.19 &  29705.34 & 0.00\%\\
\hline
\end{tabular}
\end{table}

\begin{figure}[ht]
\centering
\subfigure[Projection-based C\&CG]{
\includegraphics[width=4.2cm,height=2.1cm]{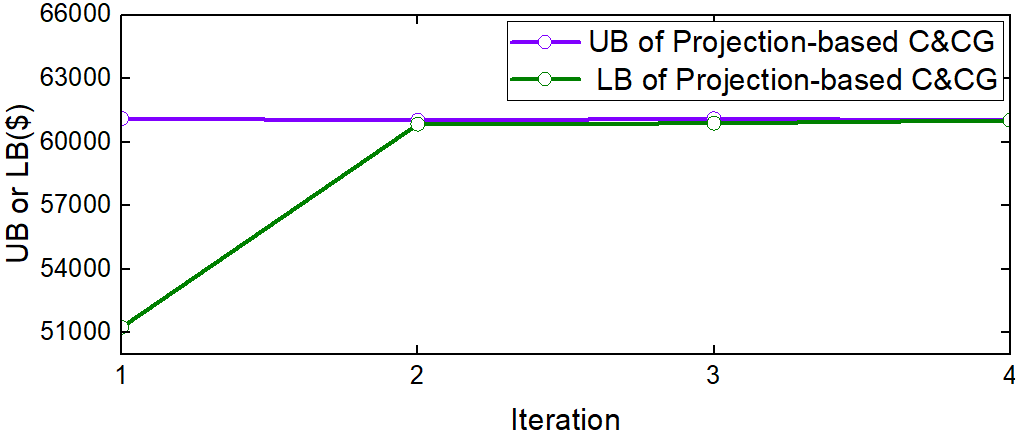}}\subfigure[Traditional C\&CG]{
\includegraphics[width=4.2cm,height=2.1cm]{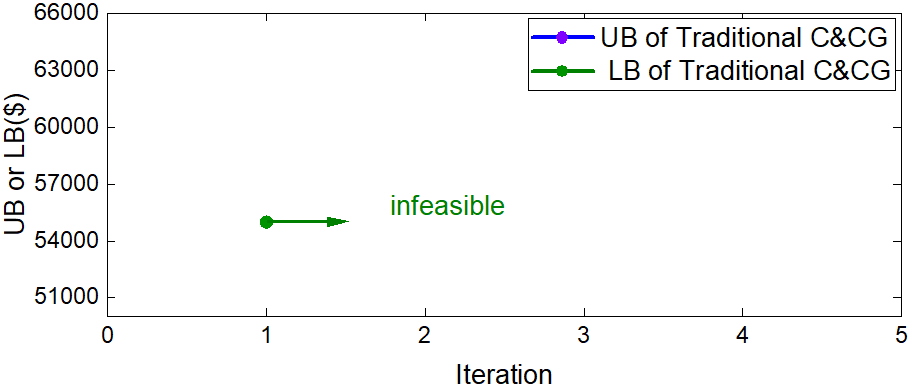}}
\caption{Iteration process of the projection-based C\&CG and traditional C\&CG algorithm.}
\label{fig:CCG}
\end{figure}

\begin{figure}[ht]
\centering
\subfigure[1st iteration]{
\includegraphics[width=4cm,height=2cm]{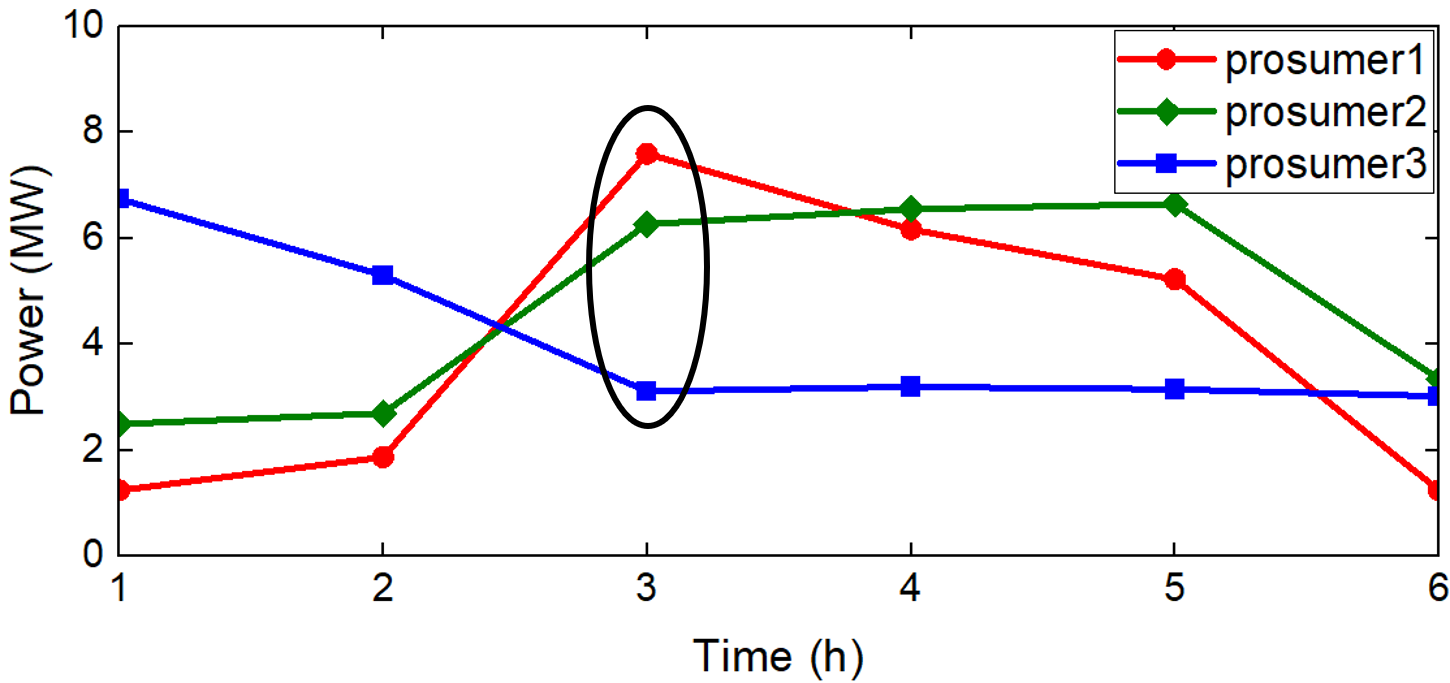}}\subfigure[2nd iteration]{
\includegraphics[width=4cm,height=2cm]{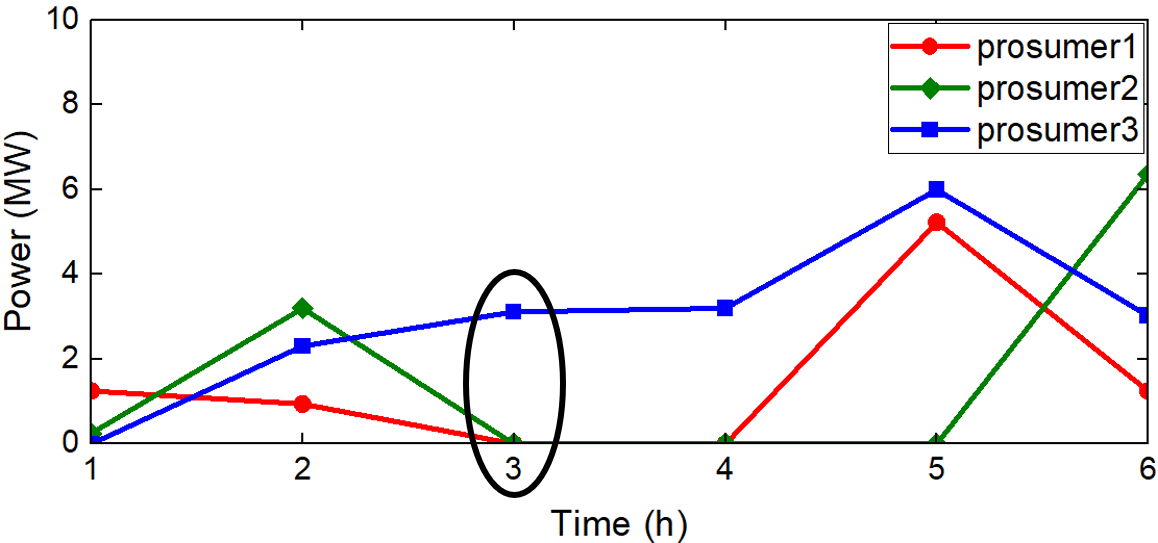}}\\
\subfigure[3rd iteration]{
\includegraphics[width=4cm,height=2cm]{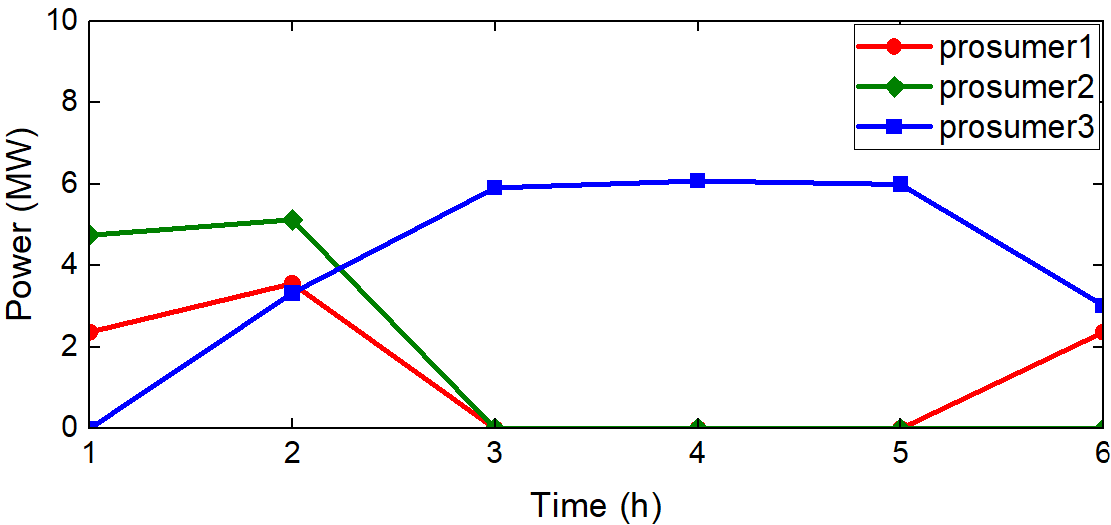}}\subfigure[4th iteration]{
\includegraphics[width=4cm,height=2cm]{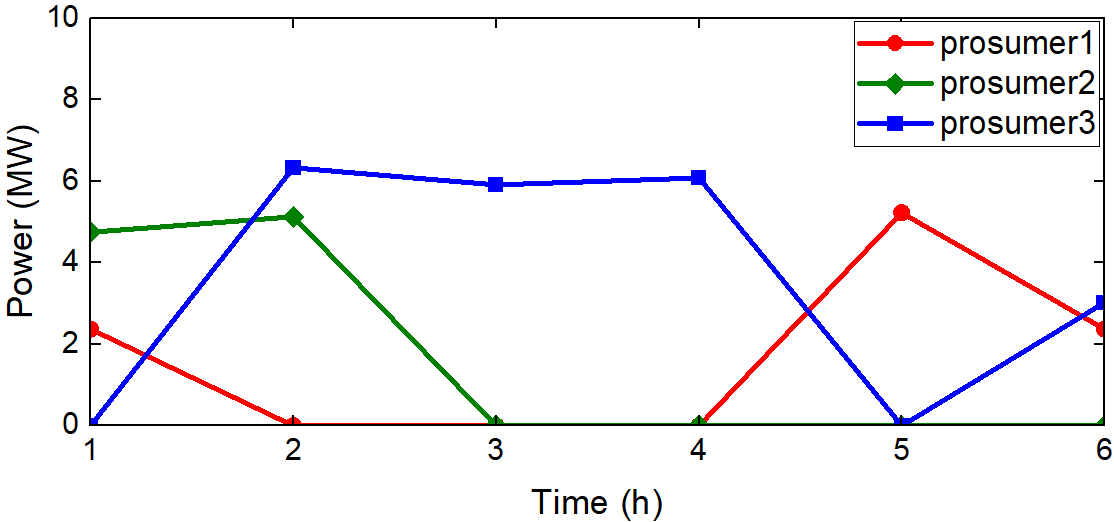}}\\
\caption{Worst-case RG scenarios during the iteration process.}
\label{fig:Scenario}
\end{figure}

To show the necessity of the projection-based C\&CG algorithm, we apply it and the traditional C\&CG algorithm to the RMD model. The changes of $UB$ and $LB$ during the iteration process are plotted in Fig. \ref{fig:CCG}. The proposed algorithm only needs 4 iterations to converge but the traditional C\&CG algorithm gets stuck at an infeasible point from the 1st iteration. The worst-case scenarios obtained in the 1st, 2nd, 3rd, and 4th iterations are shown in Fig. \ref{fig:Scenario}. 
 When the first-stage RG connection/disconnection strategies change, the generated worst-case scenario changes dramatically. For example, the worst-case scenario output for period 3 in the 1st iteration is $w=[7.6,6.3,3.1]^{T}$MW, while that in the 2nd iteration is $w=[0.0,0.0,3.1]^{T}$MW because RGs 1 and 2 are disconnected. Therefore, the previously selected scenario may fall outside the new uncertainty set and need to be projected onto the new uncertainty.

\vspace{-0.3cm}

To show the advantages of the proposed solution algorithm over the alternative DIU-based method, we compare the average number of decision variables in uncertainty sets during the iteration process, the number of iterations needed to converge, and the computation time under these two methods in Table \ref{tab:DDU-DIU}. We can find that although the two methods converge to the identical optimal value, the DIU-based method results in uncertainty sets with more decision variables, and requires more iterations and longer computation time to converge. This finding reveals the merits and necessity of the proposed method.

\begin{table}[ht]
\renewcommand{\arraystretch}{1.3}
\renewcommand{\tabcolsep}{0.5em}
\centering
\caption{Performance comparison between the proposed and the DIU-based methods}
\label{tab:DDU-DIU}
\begin{tabular}{cccccc} 
\hline 
Method  & Computation   &  Iterations &  Average number  & Converged \\
 & time(s)  &   &  of decision variables  & value  & \\
 &   &   &  in uncertainty sets  &   &  \\
\hline
Proposed  & 10.24 &  4 & 14.83 & 64139.75 \\
DIU-based & 15.51 &  7 & 18 & 64139.75\\
\hline
\end{tabular}
\end{table}

\subsection{Performance Analysis}

The impacts of two factors on the performance of the proposed method are tested, including the load flexibility and the renewable uncertainty. The robustness and scalability of the proposed model are validated.

\subsubsection{Impact of load flexibility}
First, we test how the day-ahead connection/disconnection strategy changes with a rising load flexibility reflected by the load adjustable range. The resulting underutilized RG power, the total cost \eqref{eq:robust-upper.1}, the number of iterations and the computation time are recorded in Table \ref{tab:par-impact1}. As the loads can vary within wider ranges, both the underutilized RG power and the total cost decrease. This is because higher load flexibility can help accommodate more RGs and thus reducing the underutilization cost. 

\begin{table}[ht]
\renewcommand{\arraystretch}{1.3}
\renewcommand{\tabcolsep}{0.5em}
\centering
\caption{Impact of load flexibility}
\label{tab:par-impact1}
\begin{tabular}{cccc} 
\hline 
Load range  & Underutilized RG (MW)  &  Cost (\$) & Iterations/time(s) \\
\hline
$[1.2\underline{D}_{kt},0.8\overline{D}_{kt}]$ & 42.4 &  63829 & 4/16.2 \\
$[1.1\underline{D}_{kt},0.9\overline{D}_{kt}]$ & 36.0 &  59858 & 3/13.9 \\
$[\underline{D}_{kt},\overline{D}_{kt}]$ & 36.0 &  59303 & 2/7.2 \\
$[0.9\underline{D}_{kt},1.1\overline{D}_{kt}]$ & 36.0 &  58958 & 3/9.2 \\
$[0.8\underline{D}_{kt},1.2\overline{D}_{kt}]$ & 20.0 & 49986 & 8/31.3 \\
$[0.7\underline{D}_{kt},1.3\overline{D}_{kt}]$ & 14.4 &  47093 & 5/24.8  \\
$[0.6\underline{D}_{kt},1.4\overline{D}_{kt}]$ & 14.4 &  47052 & 6/26.8 \\
$[0.5\underline{D}_{kt},1.5\overline{D}_{kt}]$ & 8.0 &  43828 & 15/375.5 \\
\hline
\end{tabular}
\end{table}

\subsubsection{Impact of renewable uncertainty}
We next change the range of RG output $[W^l, W^u]$ in the uncertainty set \eqref{eq:DDU-set}. The underutilized RG/rate, the total cost, the number of iterations and the computation time are recorded in Table \ref{tab:par-impact2}. We find that as the uncertainty range of RG output gradually enlarges, more RGs need to be disconnected to guarantee the feasibility of the real-time operation. Moreover, the growth rate of dispatch cost will become slower with the smaller variation range of renewable uncertainty. 

\begin{table}[ht]
\renewcommand{\arraystretch}{1.3}
\renewcommand{\tabcolsep}{0.5em}
\centering
\caption{Impact of renewable uncertainty}
\label{tab:par-impact2}
\begin{tabular}{cccc} 
\hline
$[W^{l},W^{u}]$  & Underutilized RG (MW)&  Total cost (\$) & Iterations\\
 & / rate(\%) &  & / time(s) \\
\hline
$[0.9,1.1]W^{e}$ & 8.0/9.0 &  52814 & 4/10.7  \\
$[0.8,1.2]W^{e}$ & 8.0/9.0 &  55648 & 2/6.7  \\
$[0.7,1.3]W^{e}$ & 24.0/26.9 & 58329 & 3/8.9  \\
$[0.6,1.4]W^{e}$ & 28.0/31.2 &  59980 & 6/14.2  \\
$[0.5,1.5]W^{e}$ & 34.0/38.1 &  61778 & 5/13.7  \\
$[0.4,1.6]W^{e}$ & 39.6/44.4 &  63153 & 6/15.2  \\
$[0.3,1.7]W^{e}$ & 46.8/52.5 &  65070 & 4/11.4  \\
\hline
\end{tabular}
\end{table}





\subsubsection{Out-of-Sample Test}

To evaluate the robustness of the strategy obtained by the proposed RMD model, out-of-sample tests are conducted. Scenarios are randomly generated from a normal distribution with the same expectation of $W^{e}$ and standard deviation of  $0.04W^{e}$, $0.06W^{e}$, $0.08W^{e}$, and $0.10W^{e}$, respectively. 500 scenarios are generated to test whether the obtained day-ahead strategy can ensure the feasibility of real-time energy sharing under the selected scenario. The percentage of infeasible scenarios and the total cost \eqref{eq:robust-eq} of three groups of uncertainty budgets are compared in Table \ref{tab:out-of-sample}. When both $B^{S}$ and $B^{T}$ equal zero, the model \eqref{eq:robust-upper} degenerates to the deterministic model without uncertainty considerations. As we can see from the table, the infeasible rate decreases with the rising $B^S$ and $B^T$. In the case with $B^S=3$ and $B^T=6$, the dispatch strategy is always feasible, even under rare extreme events. However, this is at the expense of a much higher operation cost.

\begin{table}[ht]
\renewcommand{\arraystretch}{1.3}
\renewcommand{\tabcolsep}{0.5em}
\centering
\caption{Out-of-sample test of the model with various uncertainty budgets}
\label{tab:out-of-sample}
\begin{tabular}{cccccc} 
\hline 
 Uncertainty budget &  & Infeasible & rate(\%) & & Cost(\$)\\
  &$0.04W^{e}$&$0.06W^{e}$&$0.08W^{e}$&$0.10W^{e}$& \\ 
\hline 
$B^{S}=0,B^{T}=0$ & 83.40 & 84.20 & 89.60 &88.80 & 25985\\
$B^{S}=2,B^{T}=4$ & 49.60 & 50.60 & 53.40 &53.20 & 44578\\
$B^{S}=3,B^{T}=6$ & 0.00 & 0.00 & 0.00 &0.00 & 49080\\
\hline
\end{tabular}
\end{table}

\subsubsection{Scalability}
\begin{figure}[ht]
\centering
\includegraphics[width=1\columnwidth]{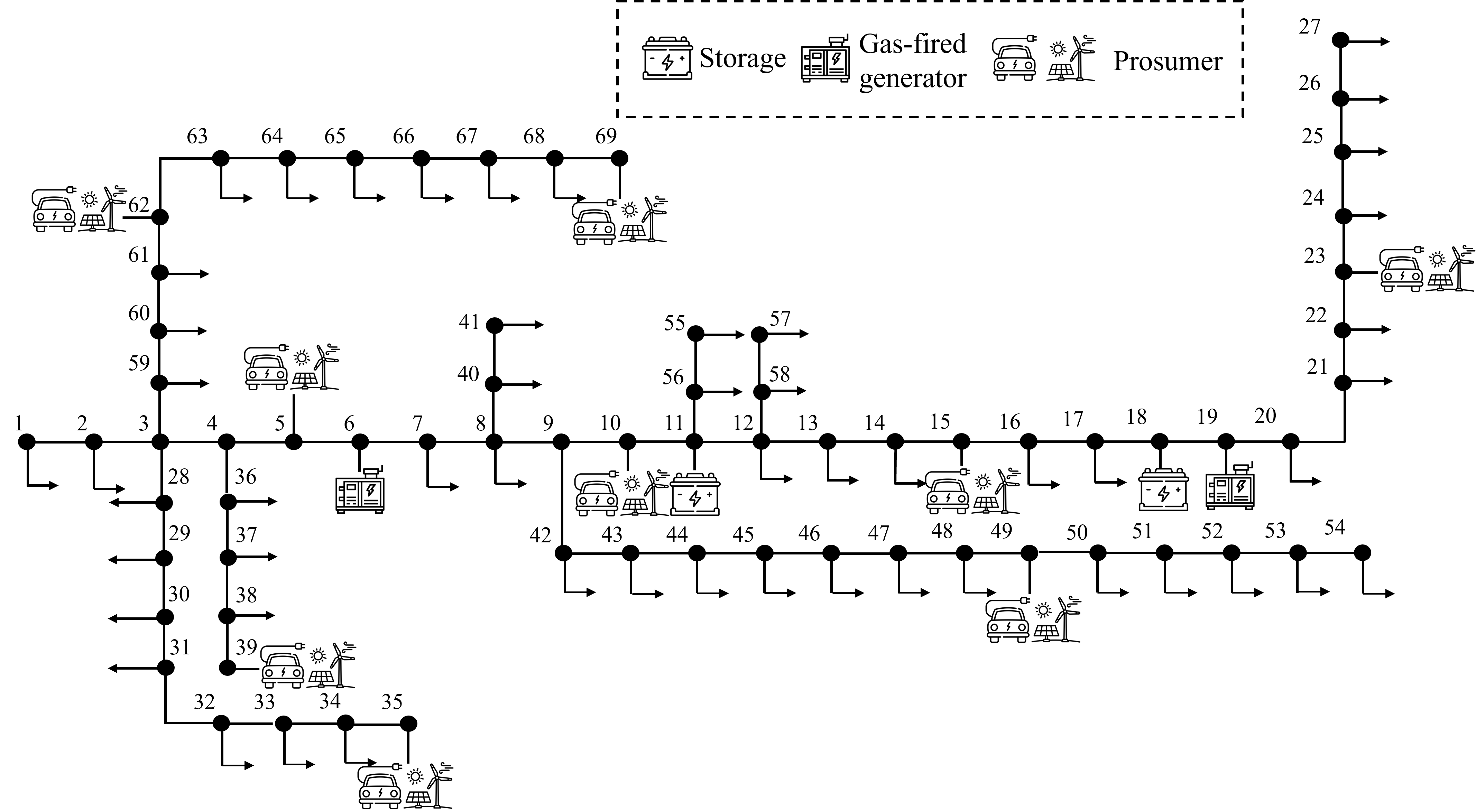}
\caption{\textcolor{black}{Topology of the modified 69-bus microgrid system with 9 prosumers.}}
\label{fig:scalability}
\end{figure}

\textcolor{black}{We use the modified 33-bus, 69-bus, 85-bus, and 123-bus microgrid systems under different settings to verify the scalability of the proposed model and algorithm. Fig.\ref{fig:scalability} shows the topology of the modified 69-bus microgrid system with 9 prosumers. The number of iterations and computational time are listed in Table \ref{tab:scabality}. Theoretically, the number of iterations to achieve convergence is less than the combination of the number of vertices of the  uncertainty set and day-ahead on/off status of renewable generators $u$. In practice, the number of iterations and computational time vary from case by case. The computational time needed is always less than 15 minutes in all test cases, which is acceptable for the day-ahead dispatch.}

\begin{table}[ht]
\renewcommand{\arraystretch}{1.3}
\renewcommand{\tabcolsep}{0.5em}
\centering
\caption{\textcolor{black}{Number of iterations / computational time}}
\label{tab:scabality}
\begin{tabular}{ccccc} 
\hline 
 \textcolor{black}{No. of prosumers} &\textcolor{black}{3} & \textcolor{black}{5} & \textcolor{black}{7} &\textcolor{black}{9}\\
\hline 
\textcolor{black}{33-bus} & \textcolor{black}{2 / 4s} & \textcolor{black}{3 / 6s} &\textcolor{black}{42 / 517s} & \textcolor{black}{45 / 895s}  \\
\textcolor{black}{69-bus} & \textcolor{black}{2 / 4s} & \textcolor{black}{3 / 10s} & \textcolor{black}{19 / 41s} & \textcolor{black}{37 / 345s}  \\
\textcolor{black}{85-bus} & \textcolor{black}{2 / 4s} & \textcolor{black}{4 / 12s} & \textcolor{black}{36 / 467s} & \textcolor{black}{26 / 271s} \\
\textcolor{black}{123-bus} & \textcolor{black}{2 / 7s} & \textcolor{black}{2 / 10s} & \textcolor{black}{62 / 821s} & \textcolor{black}{39 / 501s} \\
\hline 
\end{tabular}
\end{table}

\section{Conclusion}
 \label{sec-5}

This paper studies the robust microgrid dispatch problem with day-ahead RG connection/disconnection and real-time energy sharing between customers. A projection-based C\&CG algorithm is developed to solve the problem efficiently with convergence guarantee. The main findings are:

1) The proposed energy sharing mechanism achieves the same flexibility and efficiency as the centralized scheme.  

2) The projection-based C\&CG algorithm can address the potential failure of traditional C\&CG algorithm in dealing with endogenous uncertainty.

3) A tradeoff between security and sustainability (economy) can be observed. Basically, with more RGs being connected, the operation cost decreases while the risk due to renewable uncertainty increases.


In future research, we may consider various types of endogenous uncertainty and bounded rationality of customers.

\ifCLASSOPTIONcaptionsoff
\newpage
\fi

  
\appendices

\makeatletter
\@addtoreset{equation}{section}
\@addtoreset{theorem}{section}

\makeatother
\renewcommand{\theequation}{A.\arabic{equation}}
\renewcommand{\thetheorem}{A.\arabic{theorem}}
\section{Proof of Proposition \ref{Thm:prop-SOC}}
\label{apen-SOC}	
Let $E_{e,t-1}$ denote the SOC of energy storage units at the end of period $t-1$. It satisfies:
\begin{align}
\label{eq:apen3-1}
\underline{E}_{e,t-1}^{DA} \le E_{e,t-1} \le \overline{E}_{e,t-1}^{DA}
\end{align}

For any RG power output scenario, the real-time charging/discharging strategy $p_{et}^{c}$ and  $p_{et}^{d}$ must satisfy the constraints \eqref{eq:Fc.10}-\eqref{eq:Fc.11}. The SOC of energy storage units in period $t$ is:
\begin{align}
\label{eq:apen3-2}
E_{et} = ~& E_{e,t-1} + (p_{et}^{c}\eta_{e}^{c}-p_{et}^{d}/\eta_{e}^{d})\Delta t \nonumber\\ 
~& \le \overline{E}_{e,t-1}^{DA} +(\overline{p}_{et}^{DA,c}\eta_{e}^{c}-\underline{p}_{et}^{DA,d}/\eta_{e}^{d})\Delta t \le \overline{E}_{e,t}^{DA} \le \overline{E}_{e} \\
E_{et} = ~& E_{e,t-1} + (p_{et}^{c}\eta_{e}^{c}-p_{et}^{d}/\eta_{e}^{d})\Delta t \nonumber\\ 
~& \ge \underline{E}_{e,t-1}^{DA} +(\underline{p}_{et}^{DA,c}\eta_{e}^{c}-\overline{p}_{et}^{DA,d}/\eta_{e}^{d})\Delta t \ge \underline{E}_{e,t}^{DA} \ge \underline{E}_{e}\\
 -\Delta_{e} ~& \le \underline{E}_{eT}^{DA}-E_{e0} \le {E}_{eT}-E_{e0} \le \overline{E}_{eT}^{DA}-E_{e0} \le \Delta_{e}
\end{align}
Therefore, the SOC $E_{et}$ is always in the allowable range [$\underline{E}_{e},\overline{E}_{e}$], although SOC constraints are not included in \eqref{eq:Fc}.

\makeatother
\renewcommand{\theequation}{B.\arabic{equation}}
\renewcommand{\thetheorem}{B.\arabic{theorem}}

\section{Proof of Proposition \ref{Thm:prop-central}}
\label{apen-central}

Denote the region characterized by \eqref{eq:sharing-pro.4} as $\mathcal{D}_{kt}$. $\mathcal{D}_{kt}$ and $\mathcal{F}_{ct}$ are closed and convex. Define $\Omega_t:=\prod_{k \in \mathcal{K}} {\mathcal{D}}_{kt} \times \mathcal{F}_{ct} \times \mathbbm{R}^{(I+J)}$. Since $q_{kt}:=-a\lambda_{kt}+b_{kt}$, problem \eqref{eq:sharing-oper} can be equivalently written as
		\bsq
		\begin{align}
		\mathop{\min}_{q_{kt},\forall k \in \mathcal{K}}~ & \sum \sum \limits_{k \in \mathcal{K}} (q_{kt}-b_{kt})^2 \\
		\mbox{s.t.}~ & q_t \in \mathcal{F}_{ct}
		\end{align}
		\esq
Suppose $(d_t^*,b_t^*,q_t^*,\lambda_t^*)$ is a GNE of the energy sharing game, then according to Definition \ref{def1}, we have
\begin{align}
\label{eq:condition-2}
\forall d_t \in {\mathcal{D}}_{kt}, \forall k : ~ &  \left[
U_k(d_{kt})-U_k(d_{kt}^*)+(d_{kt}-d_{kt}^*)\lambda_{kt}^*\right] \ge 0 \nonumber\\
\forall q_t \in \mathcal{F}_{ct}:~ & \sum \limits_{k \in \mathcal{K}} (q_{kt}-q_{kt}^*)(q_{kt}^*-b_{kt}^*)\ge 0
\end{align}
		
For convex problem \eqref{eq:central-lower}, its Lagrangian function is
\begin{align}
    \mathcal{L}(d,q,\eta) & = \sum_{k \in \mathcal{K}} \sum_{t \in \mathcal{T}} U_k (d_{kt}) + \sum_{k \in \mathcal{I}} \sum_{t \in \mathcal{T}} \eta_{kt} (d_{kt} + d_{kt}^f - q_{kt}) \nonumber \\
    & + \sum_{k \in \mathcal{J}} \sum_{t \in \mathcal{T}} \eta_{kt} (d_{kt} + d_{kt}^f - w_{kt} - q_{kt}), \forall (d,q,\eta) \in \Pi_{t=1}^T \Omega_t \nonumber
\end{align}
Suppose $(\hat d,\hat q,\hat \eta)$ is a saddle point of the Lagrangian function, then $(\hat d, \hat q,\hat \eta) \in \Omega$, and it satisfies $\forall (d,q,\eta) \in \Pi_{t=1}^T \Omega_t$,
\begin{align}
    \label{eq:saddle-point}
    \mathcal{L}(d,q,\hat{\eta}) \geq \mathcal{L}(\hat{d},\hat{q},\hat{\eta}) \geq \mathcal{L}(\hat{d},\hat{q},\eta)
\end{align}
Consider three types of points $(d,q,\eta) \in \Pi_{t=1}^T \Omega_t$: 1) $d_{k't} = \hat{d}_{k't}, \forall k' \neq k, \forall t$, $q = \hat{q}$, and $\eta = \hat{\eta}$; 2) $d = \hat{d}$ and $\eta = \hat{\eta}$; 3) $d = \hat{d}$ and $q = \hat{q}$. Condition \eqref{eq:saddle-point} at these points comes down to:
		\bsq
		\label{eq:condition}
		\begin{align}
	 \left[U_{kt}(d_{kt})- U_{kt}(\hat d_{kt})+ (d_{kt}-\hat d_{kt})\hat \eta_{kt}\right] & ~ \ge 0,\forall k,\forall t \label{eq:condition.2}\\
		-\sum \nolimits_{k \in \mathcal{K}} (q_{kt}-\hat q_{kt})\hat \eta_{kt} &~\ge 0 ,\forall t\label{eq:condition.3}\\
	 \left[\sum \nolimits_{k \in \mathcal{I}}  (\eta_{kt}-\hat \eta_{kt})(\hat d_{kt} + d_{kt}^f-\hat q_{kt})\right.&\nonumber\\
		\left.+\sum \nolimits_{k \in \mathcal{J}}  (\eta_{kt}-\hat \eta_{kt})(\hat d_{kt} + d_{kt}^f-w_{kt}-\hat q_{kt})\right] &~\le 0,\forall t\label{eq:condition.4}
		\end{align}
		\esq
		
\emph{\textbf{Existence}}. Suppose $\hat d$ is the optimal solution of \eqref{eq:central-lower} and $\hat \eta$ is the corresponding dual variable. Let $d^*=\hat d$, $\lambda^*=\hat \eta$, and for all $t\in \mathcal{T}$, $q_{kt}^*=\hat d_{kt} + d_{kt}^f=\hat q_{kt}$ for all $k \in \mathcal{I}$, $q_{kt}^*=\hat d_{kt} + d_{kt}^f-w_{kt}=\hat q_{kt}$ for all $k \in \mathcal{J}$, and $b_{kt}^*=a\hat \eta_{kt}+q_{kt}^*,\forall k \in \mathcal{K}$, then it is easy to check that \eqref{eq:condition-2} is met. Thus, we have constructed a GNE $(d_t^*,b_t^*,q_t^*,\lambda_t^*)$ for each period $t$.
		
\emph{\textbf{Uniqueness}}. Given a GNE $(d_t^*,b_t^*,q_t^*,\lambda_t^*)$, when $k \in \mathcal{I}$, we have $b_{kt}^*=d_{kt}^* + d_{kt}^f+a\lambda_{kt}^*$; when $k \in \mathcal{J}$, we have $b_k^*=d_{kt}^* + d_{kt}^f-w_{kt}+a\lambda_{kt}^*$. Let $\hat d=d^*$, $\hat \eta=\lambda^*$, and $\hat q =-a\lambda^*+b^*$, then it is easy to check that $(\hat d, \hat q,\hat \eta)$ satisfies \eqref{eq:condition}, so $\hat d$ is the optimal solution of \eqref{eq:central-lower} and $\hat \eta$ is the corresponding dual variable. Since the objective function is strictly convex, and the constraint sets $\hat{D}_{kt},\forall k \in \mathcal{K},\forall t$  and $\mathcal{F}_{ct},\forall t$ are all closed convex sets, problem \eqref{eq:central-lower} has a unique solution \cite{LectureNote}, so $\hat d$ is unique.

\makeatother
\renewcommand{\theequation}{C.\arabic{equation}}
\renewcommand{\thetheorem}{C.\arabic{theorem}}
\section{Proof of Proposition \ref{Thm:prop-payment}}
\label{apen-payment}

Suppose $(d_t^*,b_t^*,q_t^*,\lambda_t^*)$ is a GNE of the energy sharing game. According to Appendix \ref{apen-central}, we know that $d^*$ satisfies \eqref{eq:condition} with $\hat \eta=\lambda^*$, $\hat q=q^*$. Therefore, if we choose $d=d^*$, $q=0$, and $\eta=\hat \eta$, \eqref{eq:condition.3} implies that $\forall t \in \mathcal{T}$
\begin{align}
    \sum \limits_{k \in \mathcal{K}} \lambda_{kt}^* q_{kt}^*=\sum \limits_{k \in \mathcal{K}}  \hat \eta_{kt} \hat q_{kt} \ge \sum \limits_{k \in \mathcal{K}} \hat \eta_{kt} q_{kt} =0,\forall t \nonumber
\end{align}

\makeatother
\renewcommand{\theequation}{D.\arabic{equation}}
\renewcommand{\thetheorem}{D.\arabic{theorem}}
\textcolor{black}{\section{Proof of Proposition \ref{Thm:prop-converge}}}
\label{apen-converge}
\textcolor{black}{Assume $(x^*,u^*)$ is an optimal solution of the proposed two-stage RO with DDU model \eqref{eq:compact} with the optimal value $\mathcal{O}^{*}$.}

\textcolor{black}{1) To begin with, we prove \textbf{Claim 1}: $LB_{j} \leq \mathcal{O}^{*} \leq UB_{j}$ for any $j$. The master problem in the $j$-th iteration is:}
\textcolor{black}{
\bsq
\label{eq:master-appendix}
\begin{align}
    \textbf{Master:}~ LB_{j} = \min_{x,u} ~ & \gamma^\top x + \beta^\top u + \tau, \\
    \mbox{s.t.}~ & Ax \ge e, u \in \{0,1\}^{|\mathcal{J}| \times |\mathcal{T}|}, \\
    ~ & H y^i \ge h- C x - T (u \circ w^i) ,\forall i \in \mathcal{S}_{w}^{j},\\
    ~ & \tau \ge \rho^\top y^i,\forall i \in \mathcal{S}_{w}^{j},
\end{align}
\esq}
\textcolor{black}{where $\mathcal{S}_{w}^{j}$ is the set of selected scenarios in the $j$-th iteration. The master problem is a relaxation of the problem \eqref{eq:compact}, which provides a lower bound for $\mathcal{O}^{*}$, i.e.,  $LB_{j} \leq \mathcal{O}^{*}$. By expanding $\mathcal{S}_{w}$ by adding identified significant scenarios from subproblems, stronger lower bounds will be gained. On the other hand, the optimality of $(x^*,u^*)$ indicates that:
\begin{align}
\label{eq:ub-appendix}
UB_{j} = \gamma^\top x^{j} + \beta^\top u^{j} +  \max_{w \in \mathcal{W}(u^{j})} \min_{y \in \Pi(x^{j},w)} \rho^\top y \ge \mathcal{O}^{*}.
\end{align}}

\textcolor{black}{2) Then we prove \textbf{Claim 2}: Let $w^{j_{1}}, w^{j_{2}}$ denote the generated significant scenarios in the $j_{1}$-th and $j_{2}$-th iterations, respectively. If the proposed projection-based C\&CG algorithm doesn't converge after $\mathcal{Z}$ iterations, then for any $j_{1}, j_{2} \leq \mathcal{Z}$, we have $w^{j_{1}} \neq w^{j_{2}}$. }

\textcolor{black}{Without loss of generality, we assume $j_{1} \textless j_{2}$. If we have $w^{j_{1}} = w^{j_{2}}$, we must have $ w^{j_{2}} = w^{j_{1}} \in  \mathcal{S}_{w}^{j_{2}}$. Therefore,}

\textcolor{black}{
\begin{align}
\label{eq:lb-ub-appendix}
\begin{aligned}
LB_{j_{2}} & = \gamma^\top x^{j_{2}} + \beta^\top u^{j_{2}} + \max_{w \in \mathcal{S}_{w}^{j_{2}}} \min_{y \in \Pi(x^{j_{2}},w)} \rho^\top y \\ 
          & \ge \gamma^\top x^{j_{2}} + \beta^\top u^{j_{2}} + \min_{y \in \Pi(x^{j_{2}},w^{j_{2}})} \rho^\top y = UB_{j_{2}}
\end{aligned}
\end{align}}

\textcolor{black} {Together with $LB_{j_{2}} \leq UB_{j_{2}}$ from \textbf{Claim 1}, we then have $LB_{j_{2}} = UB_{j_{2}}$ and thus the algorithm ends. This conclusion contradicts with the assumption that the proposed projection-based C\&CG algorithm doesn't converge after $\mathcal{Z}$ iterations.}

\textcolor{black}{3) Finally we prove \textbf{Claim 3}: The algorithm will converge to the optimal value $\mathcal{O}^*$ with an error of at most $\varepsilon$ within $\mathcal{O}(n)$ iterations. We know that the selected worst-case scenario can be achieved at a vertex of the uncertainty set in each iteration \cite{lorca2014adaptive}. According to \textbf{Claim 2}, the same vertex will not be selected twice. Hence, the projection-based C\&CG algorithm will stop within $\mathcal{O}(n)$ iterations. Suppose the algorithm converges after $Z$ iterations, according to \textbf{Claim 1}, we have $LB_{Z} \leq \mathcal{O}^{*} \leq UB_{Z}$. Together with the convergence condition $|UB_{j}-LB_{j}|\le \epsilon$, we have }
\textcolor{black}{\begin{align}
\label{eq:o-appendix}
|UB_{Z} - \mathcal{O}^{*}| \leq |UB_{Z} - LB_{Z}| \leq \epsilon
\end{align}}
\textcolor{black}{This completes the proof.}

\bibliographystyle{IEEEtran}
\bibliography{IEEEabrv,mybib}

\begin{thebibliography}{10}
\providecommand{\url}[1]{#1}
\csname url@samestyle\endcsname
\providecommand{\newblock}{\relax}
\providecommand{\bibinfo}[2]{#2}
\providecommand{\BIBentrySTDinterwordspacing}{\spaceskip=0pt\relax}
\providecommand{\BIBentryALTinterwordstretchfactor}{4}
\providecommand{\BIBentryALTinterwordspacing}{\spaceskip=\fontdimen2\font plus
\BIBentryALTinterwordstretchfactor\fontdimen3\font minus
  \fontdimen4\font\relax}
\providecommand{\BIBforeignlanguage}[2]{{%
\expandafter\ifx\csname l@#1\endcsname\relax
\typeout{** WARNING: IEEEtran.bst: No hyphenation pattern has been}%
\typeout{** loaded for the language `#1'. Using the pattern for}%
\typeout{** the default language instead.}%
\else
\language=\csname l@#1\endcsname
\fi
#2}}
\providecommand{\BIBdecl}{\relax}
\BIBdecl

\bibitem{ellabban2014renewable}
O.~Ellabban, H.~Abu-Rub, and F.~Blaabjerg, ``Renewable energy resources:
  Current status, future prospects and their enabling technology,''
  \emph{Renewable and Sustainable Energy Reviews}, vol.~39, pp. 748--764, 2014.

\bibitem{parag2016electricity}
Y.~Parag and B.~K. Sovacool, ``Electricity market design for the prosumer
  era,'' \emph{Nature energy}, vol.~1, no.~4, pp. 1--6, 2016.

\bibitem{qiu2020recourse}
H.~Qiu, H.~Long, W.~Gu, and G.~Pan, ``Recourse-cost constrained robust
  optimization for microgrid dispatch with correlated uncertainties,''
  \emph{IEEE Transactions on Industrial Electronics}, vol.~68, no.~3, pp.
  2266--2278, 2020.

\bibitem{qiu2020tri}
H.~Qiu, W.~Gu, Y.~Xu, W.~Yu, G.~Pan, and P.~Liu, ``Tri-level mixed-integer
  optimization for two-stage microgrid dispatch with multi-uncertainties,''
  \emph{IEEE Transactions on Power Systems}, vol.~35, no.~5, pp. 3636--3647,
  2020.

\bibitem{chu2021frequency}
Z.~Chu, N.~Zhang, and F.~Teng, ``Frequency-constrained resilient scheduling of
  microgrid: A distributionally robust approach,'' \emph{IEEE Transactions on
  Smart Grid}, vol.~12, no.~6, pp. 4914--4925, 2021.

\bibitem{olivella2018optimization}
P.~Olivella-Rosell, E.~Bullich-Massagu{\'e}, M.~Arag{\"u}{\'e}s-Pe{\~n}alba,
  A.~Sumper, S.~{\O}. Ottesen, J.-A. Vidal-Clos, and R.~Villaf{\'a}fila-Robles,
  ``Optimization problem for meeting distribution system operator requests in
  local flexibility markets with distributed energy resources,'' \emph{Applied
  energy}, vol. 210, pp. 881--895, 2018.

\bibitem{Website}
G.~of~South~Australia, ``Remote disconnect and reconnection of electricity
  generating plants,'' 2023,
  \url{https://www.energymining.sa.gov.au/__data/assets/pdf_file/0007/808225/2022D066388-Technical-Regulator-Guidelines-Distributed-Energy-Resources-Version-1.5-1.pdf}.

\bibitem{chen2022robust}
Y.~Chen and W.~Wei, ``Robust generation dispatch with strategic renewable power
  curtailment and decision-dependent uncertainty,'' \emph{IEEE Transactions on
  Power Systems}, vol.~38, no.~5, pp. 4640--4654, 2023.

\bibitem{wang2024application}
B.~Wang, X.~Wang, X.-P. Zhang, J.~Huang, Z.~Song, L.~Zhang, and Y.~Li,
  ``Application of decision-dependent uncertainty in power system planning and
  operation analyses: A state-of-the-art review,'' \emph{Electric Power Systems
  Research}, vol. 233, p. 110458, 2024.

\bibitem{qi2023chance}
N.~Qi, P.~Pinson, M.~R. Almassalkhi, L.~Cheng, and Y.~Zhuang,
  ``Chance-constrained generic energy storage operations under
  decision-dependent uncertainty,'' \emph{IEEE Transactions on Sustainable
  Energy}, vol.~14, no.~4, pp. 2234--2248, 2023.

\bibitem{yin2022coordinated}
W.~Yin, Y.~Li, J.~Hou, M.~Miao, and Y.~Hou, ``Coordinated planning of wind
  power generation and energy storage with decision-dependent uncertainty
  induced by spatial correlation,'' \emph{IEEE Systems Journal}, vol.~17,
  no.~2, pp. 2247--2258, 2022.

\bibitem{zeng2022two}
B.~Zeng and W.~Wang, ``Two-stage robust optimization with decision dependent
  uncertainty,'' \emph{arXiv preprint arXiv:2203.16484}, 2022.

\bibitem{zeng2013solving}
B.~Zeng and L.~Zhao, ``Solving two-stage robust optimization problems using a
  column-and-constraint generation method,'' \emph{Operations Research
  Letters}, vol.~41, no.~5, pp. 457--461, 2013.

\bibitem{nohadani2018optimization}
O.~Nohadani and K.~Sharma, ``Optimization under decision-dependent
  uncertainty,'' \emph{SIAM Journal on Optimization}, vol.~28, no.~2, pp.
  1773--1795, 2018.

\bibitem{lappas2018robust}
N.~H. Lappas and C.~E. Gounaris, ``Robust optimization for decision-making
  under endogenous uncertainty,'' \emph{Computers \& Chemical Engineering},
  vol. 111, pp. 252--266, 2018.

\bibitem{vayanos2020robust}
P.~Vayanos, A.~Georghiou, and H.~Yu, ``Robust optimization with
  decision-dependent information discovery,'' \emph{arXiv preprint
  arXiv:2004.08490}, 2020.

\bibitem{zhang2021robust}
Y.~Zhang, F.~Liu, Z.~Wang, Y.~Su, W.~Wang, and S.~Feng, ``Robust scheduling of
  virtual power plant under exogenous and endogenous uncertainties,''
  \emph{IEEE Transactions on Power Systems}, vol.~37, no.~2, pp. 1311--1325,
  2021.

\bibitem{avraamidou2020adjustable}
S.~Avraamidou and E.~N. Pistikopoulos, ``Adjustable robust optimization through
  multi-parametric programming,'' \emph{Optimization Letters}, vol.~14, no.~4,
  pp. 873--887, 2020.

\bibitem{su2022multi}
Y.~Su, F.~Liu, Z.~Wang, Y.~Zhang, B.~Li, and Y.~Chen, ``Multi-stage robust
  dispatch considering demand response under decision-dependent uncertainty,''
  \emph{IEEE Transactions on Smart Grid}, 2022.

\bibitem{tan2024robust}
T.~Tan, R.~Xie, X.~Xu, and Y.~Chen, ``A robust optimization method for power
  systems with decision-dependent uncertainty,'' \emph{Energy Conversion and
  Economics}, vol.~5, no.~3, pp. 133--145, 2024.

\bibitem{chen2023robust}
H.~Chen, X.~A. Sun, and H.~Yang, ``Robust optimization with continuous
  decision-dependent uncertainty with applications to demand response
  management,'' \emph{SIAM Journal on Optimization}, vol.~33, no.~3, pp.
  2406--2434, 2023.

\bibitem{lyu2021fully}
C.~Lyu, Y.~Jia, and Z.~Xu, ``Fully decentralized peer-to-peer energy sharing
  framework for smart buildings with local battery system and aggregated
  electric vehicles,'' \emph{Applied Energy}, vol. 299, p. 117243, 2021.

\bibitem{celik2017decentralized}
B.~Celik, R.~Roche, D.~Bouquain, and A.~Miraoui, ``Decentralized neighborhood
  energy management with coordinated smart home energy sharing,'' \emph{IEEE
  Transactions on Smart Grid}, vol.~9, no.~6, pp. 6387--6397, 2017.

\bibitem{lv2023optimal}
S.~Lv, G.~Sun, Z.~Wei, S.~Chen, and Z.~Chen, ``Optimal pricing and energy
  sharing of ev charging stations with an augmented user equilibrium model,''
  \emph{IEEE Transactions on Power Systems}, vol.~39, no.~2, pp. 4336--4349,
  2023.

\bibitem{chen2022review}
Y.~Chen and C.~Zhao, ``Review of energy sharing: Business models, mechanisms,
  and prospects,'' \emph{IET Renewable Power Generation}, vol.~16, no.~12, pp.
  2468--2480, 2022.

\bibitem{forfia2016view}
D.~Forfia, M.~Knight, and R.~Melton, ``The view from the top of the mountain:
  Building a community of practice with the gridwise transactive energy
  framework,'' \emph{IEEE Power and Energy Magazine}, vol.~14, no.~3, pp.
  25--33, 2016.

\bibitem{gjorgievski2021potential}
V.~Z. Gjorgievski, N.~Markovska, A.~Abazi, and N.~Dui{\'c}, ``The potential of
  power-to-heat demand response to improve the flexibility of the energy
  system: An empirical review,'' \emph{Renewable and Sustainable Energy
  Reviews}, vol. 138, p. 110489, 2021.

\bibitem{tushar2021peer}
W.~Tushar, C.~Yuen, T.~K. Saha, T.~Morstyn, A.~C. Chapman, M.~J.~E. Alam,
  S.~Hanif, and H.~V. Poor, ``Peer-to-peer energy systems for connected
  communities: A review of recent advances and emerging challenges,''
  \emph{Applied energy}, vol. 282, p. 116131, 2021.

\bibitem{mei2019coalitional}
J.~Mei, C.~Chen, J.~Wang, and J.~L. Kirtley, ``Coalitional game theory based
  local power exchange algorithm for networked microgrids,'' \emph{Applied
  Energy}, vol. 239, pp. 133--141, 2019.

\bibitem{yang2021optimal}
Y.~Yang, G.~Hu, and C.~J. Spanos, ``Optimal sharing and fair cost allocation of
  community energy storage,'' \emph{IEEE Transactions on Smart Grid}, vol.~12,
  no.~5, pp. 4185--4194, 2021.

\bibitem{cui2021economic}
S.~Cui, Y.-W. Wang, X.-K. Liu, Z.~Wang, and J.-W. Xiao, ``Economic storage
  sharing framework: Asymmetric bargaining-based energy cooperation,''
  \emph{IEEE Transactions on Industrial Informatics}, vol.~17, no.~11, pp.
  7489--7500, 2021.

\bibitem{xu2020data}
X.~Xu, Y.~Xu, M.-H. Wang, J.~Li, Z.~Xu, S.~Chai, and Y.~He, ``Data-driven
  game-based pricing for sharing rooftop photovoltaic generation and energy
  storage in the residential building cluster under uncertainties,'' \emph{IEEE
  Transactions on Industrial Informatics}, vol.~17, no.~7, pp. 4480--4491,
  2020.

\bibitem{chen2018analyzing}
Y.~Chen, W.~Wei, F.~Liu, Q.~Wu, and S.~Mei, ``Analyzing and validating the
  economic efficiency of managing a cluster of energy hubs in multi-carrier
  energy systems,'' \emph{Applied Energy}, vol. 230, pp. 403--416, 2018.

\bibitem{anoh2019energy}
K.~Anoh, S.~Maharjan, A.~Ikpehai, Y.~Zhang, and B.~Adebisi, ``Energy
  peer-to-peer trading in virtual microgrids in smart grids: A game-theoretic
  approach,'' \emph{IEEE Transactions on Smart Grid}, vol.~11, no.~2, pp.
  1264--1275, 2019.

\bibitem{morstyn2018bilateral}
T.~Morstyn, A.~Teytelboym, and M.~D. McCulloch, ``Bilateral contract networks
  for peer-to-peer energy trading,'' \emph{IEEE Transactions on Smart Grid},
  vol.~10, no.~2, pp. 2026--2035, 2018.

\bibitem{wang2019incentivizing}
J.~Wang, H.~Zhong, C.~Wu, E.~Du, Q.~Xia, and C.~Kang, ``Incentivizing
  distributed energy resource aggregation in energy and capacity markets: An
  energy sharing scheme and mechanism design,'' \emph{Applied Energy}, vol.
  252, p. 113471, 2019.

\bibitem{xu2019distributed}
Q.~Xu, T.~Zhao, Y.~Xu, Z.~Xu, P.~Wang, and F.~Blaabjerg, ``A distributed and
  robust energy management system for networked hybrid ac/dc microgrids,''
  \emph{IEEE Transactions on Smart Grid}, vol.~11, no.~4, pp. 3496--3508, 2019.

\bibitem{cui2018two}
S.~Cui, Y.-W. Wang, J.-W. Xiao, and N.~Liu, ``A two-stage robust energy sharing
  management for prosumer microgrid,'' \emph{IEEE Transactions on Industrial
  Informatics}, vol.~15, no.~5, pp. 2741--2752, 2018.

\bibitem{wang2022transactive}
B.~Wang, C.~Zhang, C.~Li, G.~Yang, and Z.~Y. Dong, ``Transactive energy sharing
  in a microgrid via an enhanced distributed adaptive robust optimization
  approach,'' \emph{IEEE Transactions on Smart Grid}, vol.~13, no.~3, pp.
  2279--2293, 2022.

\bibitem{ross2015multiobjective}
M.~Ross, C.~Abbey, F.~Bouffard, and G.~Jos, ``Multiobjective optimization
  dispatch for microgrids with a high penetration of renewable generation,''
  \emph{IEEE Transactions on Sustainable Energy}, vol.~6, no.~4, pp.
  1306--1314, 2015.

\bibitem{qi2016cybersecurity}
J.~Qi, A.~Hahn, X.~Lu, J.~Wang, and C.-C. Liu, ``Cybersecurity for distributed
  energy resources and smart inverters,'' \emph{IET Cyber-Physical Systems:
  Theory \& Applications}, vol.~1, no.~1, pp. 28--39, 2016.

\bibitem{hobbs2000strategic}
B.~F. Hobbs, C.~B. Metzler, and J.-S. Pang, ``Strategic gaming analysis for
  electric power systems: An mpec approach,'' \emph{IEEE transactions on power
  systems}, vol.~15, no.~2, pp. 638--645, 2000.

\bibitem{li2015demand}
N.~Li, L.~Chen, and M.~A. Dahleh, ``Demand response using linear supply
  function bidding,'' \emph{IEEE Transactions on Smart Grid}, vol.~6, no.~4,
  pp. 1827--1838, 2015.

\bibitem{bai2017distribution}
L.~Bai, J.~Wang, C.~Wang, C.~Chen, and F.~Li, ``Distribution locational
  marginal pricing (dlmp) for congestion management and voltage support,''
  \emph{IEEE Transactions on Power Systems}, vol.~33, no.~4, pp. 4061--4073,
  2017.

\bibitem{harker1991generalized}
P.~T. Harker, ``Generalized nash games and quasi-variational inequalities,''
  \emph{European journal of Operational research}, vol.~54, no.~1, pp. 81--94,
  1991.

\bibitem{huang2019adaptive}
Q.~Huang, R.~Huang, W.~Hao, J.~Tan, R.~Fan, and Z.~Huang, ``Adaptive power
  system emergency control using deep reinforcement learning,'' \emph{IEEE
  Transactions on Smart Grid}, vol.~11, no.~2, pp. 1171--1182, 2019.

\bibitem{lorca2014adaptive}
A.~Lorca and X.~A. Sun, ``Adaptive robust optimization with dynamic uncertainty
  sets for multi-period economic dispatch under significant wind,'' \emph{IEEE
  Transactions on Power Systems}, vol.~30, no.~4, pp. 1702--1713, 2014.

\bibitem{jiang2011robust}
R.~Jiang, J.~Wang, and Y.~Guan, ``Robust unit commitment with wind power and
  pumped storage hydro,'' \emph{IEEE Transactions on Power Systems}, vol.~27,
  no.~2, pp. 800--810, 2011.

\bibitem{zhao2014variable}
J.~Zhao, T.~Zheng, and E.~Litvinov, ``Variable resource dispatch through
  do-not-exceed limit,'' \emph{IEEE Transactions on Power Systems}, vol.~30,
  no.~2, pp. 820--828, 2014.

\bibitem{LectureNote}
N.~Angelia, ``Lecture notes for convex optimization,''
  \url{http://www.ifp.illinois.edu/~angelia/L5_exist_optimality.pdf}, 2008.

\end{thebibliography}

\end{document}